\documentclass[journal]{IEEEtran}

\usepackage[cmex10]{amsmath}
\usepackage{txfonts}
\usepackage{cite}
\usepackage{array}
\usepackage{graphicx}
\usepackage{paralist}
\usepackage{epstopdf}
\newtheorem{thm}{Theorem}
\newtheorem{lem}{Lemma}

\newtheorem{ass}{Assumption}

\begin{document}
\title{Efficient and Accurate Frequency Estimation of Multiple Superimposed Exponentials in Noise}

\author{Shanglin~Ye,~\IEEEmembership{Student~Member,~IEEE,}
        and~Elias~Aboutanios,~\IEEEmembership{Senior~Member,~IEEE}% <-this % stops a space
\thanks{S. Ye and E. Aboutanios are with School
of Electrical Engineering and Telecommunications, University of New South Wales, Sydney, NSW 2052, Australia. E-mail: shanglin.ye@unsw.edu.au, elias@ieee.org.}}

\maketitle

\begin{abstract}
\boldmath
The estimation of the frequencies of multiple superimposed exponentials in noise is an important research problem due to its various applications from engineering to chemistry. In this paper, we propose an efficient and accurate algorithm that estimates the frequency of each component iteratively and consecutively by combining an estimator with a leakage subtraction scheme. During the iterative process, the proposed method gradually reduces estimation error and improves the frequency estimation accuracy. We give theoretical analysis where we derive the theoretical bias and variance of the frequency estimates and discuss the convergence behaviour of the estimator. We show that the algorithm converges to the asymptotic fixed point where the estimation is asymptotically unbiased and the variance is just slightly above the Cramer-Rao lower bound. We then verify the theoretical results and estimation performance using extensive simulation. The simulation results show that the proposed algorithm is capable of obtaining more accurate estimates than state-of-art methods with only a few iterations.
\end{abstract}


\IEEEpeerreviewmaketitle

\section{Introduction}
\IEEEPARstart{E}{stimating} the frequencies of the components in sums of complex exponentials in noise is an important research problem as it arises in many applications such as radar, wireless communications and spectroscopy analysis \cite{Duda2011, Umesh1996}. The signal model given by
\begin{eqnarray} \label{sig_model} 
x(n) & = & \sum_{l = 1}^{L}  s_l(n) + w(n) \notag \\
& = & \sum_{l = 1}^{L} A_l e^{j2\pi f_l n} + w(n), \; n = 0 \ldots N-1.
\end{eqnarray}
Here $L$ is the number of components, which is assumed to be known \textit{a priori}, and $f_l \in [-0.5, 0.5]$ is the normalised frequency of the $l^{\textrm{th}}$ component. The noise terms $w(n)$ are additive Gaussian noise with zero mean and variance $\sigma^2$. The signal to noise ratio (SNR) of the $l^{\textrm{th}}$ component is $\rho_l = |A_l|^2/\sigma^2$. We note that resolving closely spaced components, as a distinct research problem itself, is not a concern of this paper.

The estimation of the frequencies of multi-tone exponentials has been the subject of intense research for many decades. The various algorithms that have been proposed to handle it, \cite{Stoica1997, Zielinski2011}, can be categorised into two types: non-parametric estimators and parametric estimators. 

Non-parametric estimators, including the traditional Capon  \cite{Capon1969}, APES \cite{Li1996} and IAA \cite{Yardibi2010} estimators exhibit high resolution, that is they can resolve closely spaced sinusoids but are computationally highly complex. Efficient implementation of these methods \cite{Glentis2008, Angel2012, Ming2011} consume $O(N^2 + K\log_2 K)$ for the computation of a length-$K$ spectrum estimate. The frequency estimation can be performed using peak picking but this requires a dense sampling of the spectrum, that is $K\gg N$, which exacerbates the computational burden.

Instead of estimating the signal spectrum, parametric estimators try to find accurate estimates of the signal parameters only. They can be further classified into time and frequency domain approaches. The time domain approaches are the more popular ones as they can achieve both high resolution as well as accurate estimation. These include subspace approaches such as Matrix Pencil (MP) \cite{Hua1990} and ESPRIT \cite{Haardt1995, Vanhuffel1994} which use the singular value decomposition (SVD) to separate the noisy signal into pure signal and noise subspaces, or iterative optimisation algorithms including IQML \cite{Bresler1986} and Weighted Least Squares (WLS) \cite{Chan2010, Sun2012} that minimise the error between the noisy and pure signals subject to different constraints. However, similar to non-parametric estimators, they suffer from high computational cost due to the SVD operation, matrix inversion and/or the eigenvalue decomposition involved, which require $O(N^3)$ for computation for large $N$. Frequency-domain parametric estimators, on the other hand, are computationally more efficient. The traditional CLEAN approach \cite{Gough1994}, which combines the maximiser of the discrete periodogram and an iterative estimation-subtraction procedure, is not desirable as the estimation error is of the same order as the reciprocal of the size of the discrete periodogram \cite{Aboutanios2011}. This makes it inaccurate when a sparse spectrum is calculated, or computationally complex for obtaining a dense spectrum. 

A number of efficient algorithms have been proposed in \cite{Aboutanios2005, Quinn1994, Yang2011} to refine the maximiser of the $N$-point periodogram by interpolation on Fourier coefficients. But, as these are developed for single-tone signals, they perform poorly in the multiple component case due to the bias resulting from the interference of the components with one another. Much work, \cite{Duda2011, Duda2014, Diao2014}, has been done to reduce the effect of the interference by applying the interpolators after pre-multiplying the signal by a time domain window. However, non-rectangular windows lead to deterioration of estimation accuracy. In this paper, we overcome the aforementioned limitations by proposing a novel parametric estimation algorithm that operates in the frequency domain and achieves excellent performance. The new approach is more computationally efficient than the non-parametric and time domain parametric estimators, yet it outperforms them in terms of estimation accuracy. 

The rest of the paper is organised as follows. In Section \ref{sec:Method}, we present the novel frequency estimation algorithm. We analyse the algorithm in Section \ref{sec:theo} and give its theoretical performance. In Section \ref{sec:Simu}, we show simulation results by comparing the proposed algorithm with state-of-art parametric estimators and the Cramer-Rao Lower Bound (CRLB). Finally, some conclusion are drawn in Section \ref{sec:Conc}.

\section{The Proposed Method} \label{sec:Method}
Let $\hat{\lambda}$ denote the estimate of the parameter $\lambda$. The A\&M estimator of \cite{Aboutanios2005, Aboutanios2011} is a powerful and efficient algorithm for the estimation of the frequency of a single-tone signal. It uses a two stage strategy, obtaining first a coarse estimate from the maximum of the $N$-point periodogram
\begin{eqnarray} \label{eq:max_bin}
\hat{m} & = & \arg \max_{k}\left|X(k)\right|^2 \notag \\
& = & \arg \max_{k}\left|\sum_{n = 0}^{N-1} x(n)e^{-j\frac{2\pi}{N}kn}\right|^2.
\end{eqnarray}
In the noiseless case, we have $\lim_{N\rightarrow \infty} \hat{m} = m,$ a.s. \cite{Aboutanios2005}. When $\hat{m}$ is the true maximum bin, the frequency of the signal can be written as
\begin{equation} \label{eq:f}
f = \frac{\hat{m} + \delta}{N},
\end{equation}
where $\delta \in [-0.5,0.5]$ is the frequency residual.  The A\&M algorithm then refines the coarse estimate by interpolating on Fourier coefficients to obtain an estimate for the residual $\delta$. 

Let $X_{\pm 0.5}$ be the Fourier Coefficients at locations $\hat m \pm 0.5$. In the noiseless case, these are given by 
\begin{eqnarray} \label{eq:unbiased_X}
X_{\pm 0.5} & = & \frac{1}{N}\displaystyle \sum\limits_{n=0}^{N-1} x(n) e^{- j\frac{2\pi}{N}\left(\hat m \pm 0.5\right)n}\notag \\
& = & \displaystyle\frac{A}{N} \frac{ 1 + e^{j2\pi\delta} }{1- e^{ j\frac{2\pi}{N}\left(\delta \mp 0.5\right)}}.
\end{eqnarray}
Putting $z^{-1} = e^{-j2\pi\frac{\delta}{N}}$, an estimate of which is constructed as
\begin{eqnarray} \label{eq:gest}
\hat{z}^{-1} = \cos\left(\frac{\pi}{N}\right)-2jh\sin\left(\frac{\pi}{N}\right),
\end{eqnarray}
where $h$ is the interpolation function
\begin{equation}
h = \frac{1}{2} \frac{X_{0.5}+X_{-0.5}}{X_{0.5}-X_{-0.5}}.
\end{equation}
From $\hat{z}^{-1}$, estimates of $\delta$ and consequently of the frequency $f$ become
\begin{equation} \label{eq:est}
\hat{\delta} = \displaystyle -\frac{N}{2\pi}\Im\left\{\ln\hat{z}^{-1}\right\},\text{ and } \hat f = \frac{\hat m+\hat \delta}{N}.
\end{equation}
Here $\Im\{\bullet \}$ denotes the imaginary part of $\bullet$. The key to the A\&M algorithm compared to other interpolators like those of \cite{Quinn2001} is that it can be implemented iteratively in order to improve the estimation accuracy, \cite{Aboutanios2005}. In each iteration the estimator removes the previous estimate of the residual before re-calculating Fourier coefficients and re-interpolating. It was shown in \cite{Aboutanios2005} that two iterations are sufficient for the estimator to obtain asymptotically unbiased frequency estimate with the variance only $1.0147$ times the CRLB.

Now turning to the multi-tone case, that is $L \geq 2$, let $\{\hat{m}_l\}_{l = 1}^{L}$ be the estimates of the maximum bins. Also let $\hat{\delta}_p$ be the estimates of the residuals from the previous iteration. The Fourier coefficients of the $p^{\textrm{th}}$ component at locations $(\hat{m}_p + \hat{\delta}_p \pm 0.5)$ are
\begin{eqnarray} 
X_{p,\pm 0.5} & = &  \frac{1}{N}\displaystyle \sum\limits_{n=0}^{N-1} x(n) e^{- j\frac{2\pi}{N}\left(\hat{m}_p + \hat{\delta}_p \pm 0.5\right)n}\label{eq:biased_coef1} \\
& = &  S_{p,\pm 0.5} + \sum_{l = 1, l\neq p}^{L} S_{l,\pm 0.5} + W_{p,\pm 0.5}\label{eq:biased_coef},
\end{eqnarray}
where $S_{p,\pm 0.5}$ are Fourier coefficients for a single exponential $s_p(n)$ as per Eq. (\ref{eq:unbiased_X}). $W_{p, \pm 0.5}$ are the corresponding noise terms at the interpolation locations. $S_{l,\pm 0.5}$ is the leakage term introduced by the $l^{\textrm{th}}$ component,
\begin{eqnarray} \label{eq:leakage}
S_{l,\pm 0.5}& = & \frac{1}{N}\sum_{n = 0}^{N-1} s_l(n) e^{- j\frac{2\pi}{N}(\hat{m}_p +\hat{\delta}_p \pm 0.5)n}\notag \\
& = & \displaystyle \frac{A_{l}}{N} \frac{ 1 + e^{j2\pi (\delta_l - \hat{\delta}_p)}}{1- e^{j\frac{2\pi}{N}\left(M_{l,p} + \delta_l - \hat{\delta}_p \mp 0.5\right)}},
\end{eqnarray}
where
\begin{equation} \label{eq:Delta}
M_{l,p} = \hat{m}_l - \hat{m}_p.
\end{equation}

As proposed in  \cite{Duda2011, Belega2014},  the leakage terms can be attenuated by applying a window to the signal. Although this reduces the bias it also leads to a broadening of the main lobe and comes at the cost of an increase in the estimation variance. We, on the other hand, address this problem by estimating the leakage terms in Eq. (\ref{eq:leakage}) and removing them in order to obtain the expected coefficients of a single exponential. We then apply the A\&M estimator to estimate the frequency. It is clear from Eq. (\ref{eq:leakage}) that this necessitates the parameters $\delta_l$ and $A_l$ be known or at least estimates for them should be available. In what follows, we construct a procedure to achieve this. 

Let us start by assuming that we have estimates $\left\{\hat{\delta}_l, \hat{A}_l \right\}_{l = 1, l \neq p}^{L}$. Then the total leakage term can be obtained as 
\begin{eqnarray} \label{eq:Total_Leakage_Term}
\sum_{l = 1, l\neq p}^L \hat{S}_{l,\pm 0.5}&=&\sum_{l = 1, l\neq p}^L \frac{\hat A_{l}}{N} \frac{ 1 + e^{ j2\pi (\hat{\delta}_l - \hat{\delta}_p)} }{1- e^{j\frac{2\pi}{N}\left(M_{l,p} + \hat{\delta}_l - \hat{\delta}_p \mp 0.5\right)}}.
\end{eqnarray}
Subtracting the estimated total leakage from the Fourier Coefficient (shown in Eq. (\ref{eq:biased_coef})) yields
\begin{eqnarray} \label{eq:unbiased_coef_est}
\hat{S}_{p,\pm 0.5} & = & X_{p,\pm 0.5} - \sum_{l = 1, l\neq p}^L \hat{S}_{l,\pm 0.5}\notag \\
& = & X_{p,\pm 0.5} - \sum_{l = 1, l\neq p}^L  \frac{\hat A_{l}}{N}  \frac{ 1 + e^{ j2\pi (\hat{\delta}_l - \hat{\delta}_p)} }{1- e^{j\frac{2\pi}{N}\left(M_{l,p} + \hat{\delta}_l - \hat{\delta}_p \mp 0.5\right)}}.
\end{eqnarray}

Accurately estimating and subtracting the leakage terms would lead to a reduction in the bias in the estimates of $\delta_p$ and $A_p$. Therefore, we propose an iterative procedure where in each iteration we estimate the parameters, $\delta$ and $A$, of each component using the previous estimates of all components. We start by initialising all of the parameter estimates to 0. To elucidate the procedure, consider the estimation of the $p^{\textrm{th}}$ component during the $i^{\textrm{th}}$ iteration. Given the set of estimates $\left\{\hat{\delta}_l^{(i-1)},\hat{A}_l^{(i-1)}\right\}_{l=1, l \neq p}^L$, we calculate the total leakage term at the locations $\hat m_p+\hat{\delta}_p^{(i-1)}\pm 0.5$ according to Eq. (\ref{eq:Total_Leakage_Term}). We then obtain the ``leakage-free'' Fourier coefficients using Eq. (\ref{eq:unbiased_coef_est}) and apply the A\&M algorithm to get a new estimate of the residual $\hat{\delta}^{(i)}_p$. The estimate of the complex amplitude, on the other hand, is obtained by subtracting the sum of the leakage terms from the amplitude at the true frequency:
\begin{eqnarray}\label{eq:est_A}
\hat{A}_p & = & \frac{1}{N}\left(\sum_{n=0}^{N-1} x(n) e^{- j2\pi\hat{f}_pn} - \sum_{l = 1, l \neq p}^{L} \hat{S}_{l, \hat{f}_l}\right) \notag \\
& = &\frac{1}{N}\sum_{n=0}^{N-1} x(n) e^{- j2\pi\hat{f}_pn} - \sum_{l = 1, l \neq p}^{L}\frac{\hat{A}_{l}}{N}\frac{ 1 - e^{j2\pi N(\hat{f}_{l}-\hat{f}_p)}}{1- e^{ j2\pi(\hat{f}_{l} -\hat{f}_p)}}.
\end{eqnarray}

As we show in the following section, as the algorithm is iterated, the error between $\hat{S}_{p, \pm 0.5}$ and $S_{p, \pm 0.5}$ is consistently reduced since  the leakage terms in Eq. (\ref{eq:leakage}) is better estimated. As the number of iterations increases, the estimation variances approach the CRLB of the multiple component case.

The estimation procedure of the proposed algorithm is summarised in Table \ref{tab:TIME}. Its overall computational complexity is $O(LN\log_2 N)$. Asymptotically, this is of the same order as the FFT operation. It is therefore more efficient than the non-parametric estimators and the time-domain parametric estimators, especially when $N$ becomes so large that the SVD, matrix inversion and eigenvalue decomposition operations utilised in those methods become unimplementable.

\begin{table}[h!]
	\caption{Estimation Procedure of the Iterative Multiple Component Estimator}
	\label{tab:TIME}
	\centering
	\begin{tabular}{ll}
     \hline
     \hline
		\textbf{1.} & Initialise $\{\hat{f}_l\}_{l=1}^L = \{\hat{\delta}_l\}_{l=1}^L = \{\hat{A}_l\}_{l=1}^L = 0$;\\
		\textbf{2.} & For $q = 1$ to $Q$ iterations do:\\
		& \quad For $p = 1$ to $L$, do:\\
		& \qquad (1) If $q=1$, find the maximum bin $\hat{m}_p$;\\
		& \qquad (2) Use Eqs. (\ref{eq:biased_coef1}) and (\ref{eq:unbiased_coef_est}) to obtain the ``leakage-free''\\
		&\qquad\qquad Fourier Coefficients $\hat{S}_{p,\pm 0.5}$;\\
		& \qquad (3) Apply the A\&M estimator (Eqs. (\ref{eq:gest})-(\ref{eq:est})) using $\hat{S}_{p,\pm 0.5}$ \\
		& \qquad \qquad to get $\hat{\delta}_p$;\\
		& \qquad (4) Update $\hat{A}_p$ using Eq. (\ref{eq:est_A});\\
		\textbf{3.}& Finally obtain $\displaystyle \hat{f}_l = \frac{\hat{m}_l + \hat{\delta}_l}{N}$, $l=1\dots L$.\\
       \hline
       \hline
	\end{tabular}
\end{table}

%%%%%%%%%%%%%%%%%%%%%%%%%%%%%%%%%%%%%%%%%%%%%%%%%%%%%%%%%%%%%%%%%%%%%%%%%%%%%%

\section{Analysis} \label{sec:theo}
In this section, we present the theoretical analysis of the proposed algorithm. We proceed to derive the theoretical bias and variance, then discuss the convergence properties. Although the noise terms in Eq. (\ref{sig_model}) are assumed to be additive Gaussian, the following analysis works under more general assumption on the noise terms established in \cite{An1983}. Under these assumptions, the Fourier coefficients of the noise converge in distribution, satisfying
\begin{equation}\label{eq:Noise_ass}
\limsup\limits_{N\rightarrow \infty} \sup\limits_{k} \frac{|W(k)|^2}{N^{-1} \ln N} \leq 2\pi f_x(k),
\end{equation}
where $f_x(k)$ is the spectrum density function of the noise, \cite{An1983}. Thus $W(k) \sim \left(N^{-\frac{1}{2}}\sqrt{\ln N}\right)$ almost surely (a.s.). To proceed, we make the following assumption on the frequency separations
\begin{ass} 
For $L \geq 2$, we have 
$$\Delta = \min\limits_{p, l = 1\ldots L; p\neq l}\left|f_p - f_l\right| \sim O(1).$$
\end{ass}
This assumption implies that the minimum frequency separation is independent of $N$.

\subsection{Theoretical Bias and Variance}
We first carry out analysis for $L = 2$, and then generalise the results to $L\geq2$. 

Let $\nu_1$ and $\nu_2$ be the estimates of $\delta_1$ and $\delta_2$ respectively, and $M = M_{2,1} = \hat{m}_{2} - \hat{m}_{1}$. Putting $V_{\pm} = S_{2,\pm 0.5} - \hat{S}_{2,\pm 0.5}$ and replacing  the subscripts $\lbrace 1, \nu_1 \pm 0.5 \rbrace$ by $\lbrace \pm \rbrace$, the ``leakage free" Fourier coefficients of the $1^{\textrm{st}}$ component can be expressed as
\begin{eqnarray}  
\hat{S}_{\pm}  & = &  S_{\pm} + V_{\pm} + W_{\pm} \label{eq:interp_coef}\\
& = &  \frac{a_1}{N\left( 1 - e^{j\frac{2\pi}{N}(\delta_1-\nu_1)} \right)} + V_{\pm} + W_{\pm} \label{eq:interp_coef_approx2}
\end{eqnarray} 
where
\begin{eqnarray} \label{eq:a}
a_1 & = & A_1\left(1 + e^{j2\pi(\delta_1-\nu_1)}\right) \notag \\
& = & 2A_1\cos(\pi (\delta_1-\nu_1)) e^{j\pi (\delta_1-\nu_1)}.
\end{eqnarray}

The interpolation function of the 1$^{\textrm{st}}$ component is
\begin{eqnarray} \label{eq:h1}
\hat{h}_1 & = & \frac{1}{2}\frac{\hat{S}_{+} + \hat{S}_-}{\hat{S}_{+} - \hat{S}_{-}} \notag\\
& = & \frac{1}{2}\frac{\frac{S_{+} + S_-}{S_+ - S_-} + \frac{V_+ + V_- + W_{+}+W_{-}}{S_+ - S_-} }{1 +  \frac{V_+ - V_- + W_{+} - W_{-}}{S_+ - S_-}} \notag \\
&=& \frac{h_1 +j\frac{N\lambda_1}{4a_1z_1\sin\left(\frac{\pi}{N}\right)}\left(V_+ + V_- + W_{+}+W_{-}\right)}
	{1+j\frac{N\lambda_1 }{2a_1z_1\sin\left(\frac{\pi}{N}\right)}\left(V_+ - V_- + W_{+}-W_{-}\right)},
\end{eqnarray}
where
\begin{equation}
h_1 = \frac{1}{2} \frac{S_{0.5}+S_{-0.5}}{S_{0.5}-S_{-0.5}},
\end{equation}
\begin{equation}
z_1 = e^{j\frac{2\pi}{N}(\delta_1 - \nu_1)},
\end{equation}
and
\begin{eqnarray}
\lambda_1 = \left(1-z_1e^{-j\frac{\pi}{N}}\right)\left(1-z_1 e^{j\frac{\pi}{N}} \right) \label{eq:lambda}.
\end{eqnarray}

Now
\begin{eqnarray}  \label{eq:V}
V_{\pm}  =  \frac{A_2}{N} \frac{1+e^{j2\pi(\delta_2-\nu_1)}}{1-e^{j\frac{2\pi}{N}\left(M + \delta_2 - \nu_1 \mp 0.5\right)}}  - \frac{\hat{A}_2}{N} \frac{1+e^{j2\pi(\nu_2-\nu_1)}}{1-e^{j\frac{2\pi}{N}\left(M + \nu_2 - \nu_1 \mp 0.5\right)}},%\notag \\
\end{eqnarray}
where $\hat{A}_2$ can be expanded as
\begin{eqnarray}\label{eq:A2}
\hat{A}_2  & = & \frac{1}{N}\sum_{n=0}^{N-1} [s_1(n) + s_2(n) + w(n) - \hat{s}_1(n)] e^{- j2\pi\hat{f}_2n}  \notag \\
& = & A_2 \frac{1 - e^{j2\pi(\delta_2 - \nu_2)}}{N\left(1 - e^{j\frac{2\pi}{N}(\delta_2 - \nu_2)}\right)} + A_1\frac{1 - e^{j2\pi(\delta_1 - \nu_2)}}{N\left(1 - e^{j\frac{2\pi}{N}(\delta_1 - \nu_2 - M)}\right)} \notag \\
&& \; - \hat{A}_1 \frac{1 - e^{j2\pi(\nu_1 - \nu_2)}}{N\left(1 - e^{j\frac{2\pi}{N}(\nu_1 - \nu_2 - M)}\right)} + W_2.
\end{eqnarray}
In Eq. (\ref{eq:A2}), $W_2$ is the noise term at $\hat{f_2} = (\hat{m}_2 + \nu_2)/N$. Under \textit{Assumption 1}, we have that $|f_2 - f_1| = |(M + \delta_2 - \delta_1)/N| \sim  O(1)$, so $M \sim O(N)$ and  $M/N \sim O(1)$. As a result, the second and third terms of Eq. (\ref{eq:A2}) are $O(N^{-1})$, and
\begin{eqnarray} \label{eq:A2_fin}
\hat{A}_2 & = & A_2 \frac{1 - e^{j2\pi(\delta_2 - \nu_2)}}{N\left(1 - e^{j\frac{2\pi}{N}(\delta_2 - \nu_2)}\right)} + W_2  + O(N^{-1}).
\end{eqnarray}

Therefore
\begin{eqnarray}  \label{eq:V_final} 
V_\pm & = &  \frac{A_2}{N}\beta_{\pm} - \frac{1+e^{j2\pi (\nu_2-\nu_1)}}{N\left(1-e^{j\frac{2\pi}{N}(M+\nu_2 - \nu_1 \mp 0.5)}\right)} ( W_2 + O(N^{-1}))\notag \\
& = &   \frac{A_2}{N} \beta_{\pm} - W_{2\pm} + O(N^{-2}),
\end{eqnarray}
where
\begin{eqnarray} \label{eq:beta1}
\beta_{\pm} & = & \frac{1+e^{j2\pi(\delta_2-\nu_1)}}{1-e^{j\frac{2\pi}{N}\left(M + \delta_2 - \nu_1 \mp 0.5\right)}} \notag \\
&& \; - \frac{\left(1 - e^{j2\pi(\delta_2 - \nu_2)}\right)\left(1+e^{j2\pi(\nu_2-\nu_1)}\right)}{N\left(1 - e^{j\frac{2\pi}{N}(\delta_2 - \nu_2)}\right)\left(1-e^{j\frac{2\pi}{N}\left(M + \nu_2 - \nu_1 \mp 0.5\right)}\right)},
\end{eqnarray}
\begin{equation} \label{eq:W2}
W_{2\pm} = \eta_\pm W_2,
\end{equation}
and
\begin{equation} \label{eq:eta}
\eta_\pm = \frac{1+e^{j2\pi (\nu_2-\nu_1)}}{N\left(1-e^{j\frac{2\pi}{N}(M+\nu_2 - \nu_1 \mp 0.5)}\right)}.
\end{equation}

Now we have the following lemma:
\begin{lem}
For $L=2$, if \textit{Assumption 1} holds, $\hat{f}_1$ asymptotically converges to a normal:
\begin{equation}
\hat{f}_1 \rightarrow \mathcal{N}(f_1 + \mu_{f_1}, \sigma^2_{f_1}).
\end{equation}
The mean $\mu_{f_1}$ and variance $\sigma^2_{f_1}$ of the distribution are respectively given by
\begin{eqnarray} \label{eq:E}
\mu_{f_1}  & = &  \frac{\pi[(\delta_1 - \nu_1)^2 - 0.25]}{2 A_1 \cos(\pi(\delta_1-\nu_1))} \notag \\
&& \; \times \left\{ [1-2(\delta_1-\nu_1)]\Im \left\{e^{-j2\pi(\delta_1 - \nu_1)}\beta_+\right\}\right. \notag \\
&& \qquad \left. + [1+2(\delta_1-\nu_1)]\Im \left\{e^{-j2\pi(\delta_1 - \nu_1)}\beta_-] \right\} \right\},
\end{eqnarray}
and 
\begin{eqnarray}\label{eq:Var}
\sigma^2_{f_1} =  \frac{\pi^2  [(\delta_1-\nu_1)^2-0.25]^2}{4\rho_1 N^3\cos^2(\pi (\delta_1-\nu_1))} [1+4(\delta_1-\nu_1)^2].
\end{eqnarray}
Furthermore, $\mu_{f_1}|_{\delta_1 = \nu_1, \delta_2 = \nu_2} = 0$, meaning the estimator is unbiased at the true frequencies.

\textit{Proof:} See Appendix A.
\end{lem}

Based on \textit{Lemma 1}, now we have the following theorem:
\begin{thm}
If \textit{Assumption 1} holds, the frequency estimates given by the proposed algorithm are asymptotically statistically independent and converge in distribution to the standard normal. that is
\begin{equation}
\forall p = 1\ldots L, \quad \hat{f}_p \rightarrow \mathcal{N}\left(f_p + \mu_{f_p},  \sigma^2_{{f}_p}\right),
\end{equation}
where  $\mu_{f_p}$ and $\sigma^2_{{f}_p}$ are respectively the mean and variance of the asymptotic distribution of $\hat{f}_p$.
\end{thm}

From Eq. (\ref{eq:E}) it is straightforward to see that the mean of the asymptotic distribution of $\hat{f}_p$ is given by
\begin{eqnarray} \label{eq:E_gen}
\mu_{f_p} & = &  \frac{\pi[(\delta_p - \nu_p)^2 - 0.25]}{2 A_p \cos(\pi(\delta_p-\nu_p))} \notag \\
&& \: \times \left\{ [1-2(\delta_p-\nu_p)]\Im \left\{e^{-j2\pi(\delta_p - \nu_p)}  \sum_{l = 1, l \neq p}^L \beta_{+,l}\right\}\right. \notag \\
&& \; \left. + [1+2(\delta_p-\nu_p)]\Im \left\{e^{-j2\pi(\delta_p - \nu_p)}\sum_{l = 1, l \neq p}^L \beta_{-,l} \right\} \right\},
\end{eqnarray}
where $\beta_{\pm, l}$ are obtained by substituting the appropriate parameters into Eq. (\ref{eq:beta1}). This result implies $\mu_{f_p} |_{\delta_l = \nu_l (l = 1\ldots L)} = 0$. 

The asymptotic variance, on the other hand, is given by:
\begin{eqnarray} \label{eq:var_gen}
\sigma^2_{{f}_p} =  \frac{\pi^2  [(\delta_p-\nu_p)^2-0.25]^2}{4\rho_p N^3\cos^2(\pi (\delta_p-\nu_p))} [1+4(\delta_p-\nu_p)^2].
\end{eqnarray}

\subsection{Convergence}
We now turn to the convergence of the proposed estimator. We have the following theorem:
\begin{thm}
If \textit{Assumption 1} holds, the proposed algorithm converges to the fixed point with the following properties:
\begin{enumerate}
\item The $L$-dimensional asymptotic fixed point is at the true frequency residuals $(\delta_1,\ldots, \delta_L)$;
\item The convergence rate is given by 
\begin{eqnarray} \label{eq:Conv_rate}
r_L \leq \sqrt{L} O\left( \max_{ l,p = 1\ldots L, l \neq p }\left\{\sqrt{\Gamma_{l,p}},  \sqrt{\frac{\ln N}{N}} \right\}\right), \quad \textrm{a.s.}
\end{eqnarray}
where
\begin{eqnarray}\label{eq:Gamma2}
\Gamma_{l,p} =  \frac{|A_l|^2}{\langle N(f_l - f_p) \rangle ^2} \sim |A_l|^2 O\left(N^{-2} \right),
\end{eqnarray}
and $\langle \bullet \rangle$ finds the nearest integer of $\bullet$.
\item Asymptotically, the estimator is unbiased and the theoretical variance of the distribution of $\hat{f}_p$ is 
\begin{eqnarray} \label{eq:varf_fin}
\sigma^2_{\hat{f}_p} & = &  \frac{\pi^2}{64\rho_p N^3},
\end{eqnarray}
which is approximately $1.0147$ times the asymptotic CRLB (ACRLB).
\end{enumerate}

\textit{Proof:} See Appendix B.
\end{thm}

Looking at the convergence rate, we see from Appendix B that the first term in the maximum in Eq. (\ref{eq:Conv_rate}) is due to the leakage terms $V_\pm$ while the second term results from the noise. Asymptotically, $r_L = O(N^{-\frac{1}{2}}\sqrt{\ln N})$ a.s. and the convergence is dictated by the noise, which is similar to the single component case. Thus two iterations are sufficient for the estimation error of residuals to become lower order than the ACRLB \cite{Yao1995}.

For finite $N$, on the other hand, the convergence behaviour is dictated by $\Gamma_{l,p}$. Consider the equality implied by Eq. (\ref{eq:Conv_rate}):
\begin{eqnarray}
 \max\limits_{l,p; l \neq p}\left\{\Gamma_{l,p}\right\} = \frac{\ln N}{N}.
\end{eqnarray}
As $N$ increases such that $\max_{l,p; l \neq p}\left\{\Gamma_{l,p}\right\} < N^{-1}\ln N$, we have $r_L \leq \sqrt{L} O(N^{-\frac{1}{2}}\sqrt{\ln N})$ a.s.. For $L \ll N$, the convergence rate is still dictated by the noise, with two iterations sufficing. When $\max_{l,p; l \neq p}\left\{\Gamma_{l,p}\right\} \geq N^{-1}\ln N$, we have $r_L \leq \sqrt{L} O\left(\max_{ l,p; l \neq p} \left\{\sqrt{\Gamma_{l,p}}\right\}\right)$. To conclude, we have the following theorem:
\begin{thm}
If \textit{Assumption 1} holds, the convergence rate $r_L$ of the proposed algorithm on a signal with $L$ components is given for the following cases:
\begin{enumerate}
\item Asymptotically:
$$r_L = O\left(\sqrt{\frac{\ln N}{N}}\right) \quad \textrm{a.s.},$$
and the algorithm converges in two iterations;
\item For finite $N$:
\begin{enumerate}
\item When $\max\limits_{l,p; l \neq p}\left\{\frac{|A_l|^2}{ \langle N(f_l - f_p)\rangle^2}\right\} < \frac{\ln N}{N}$:
$$r_L \leq \sqrt{L} O\left(\sqrt{\frac{\ln N}{N}}\right) \quad \textrm{a.s.},$$
and the algorithm converges in two iterations provided $L \ll N$;
\item When $\max\limits_{l,p; l \neq p}\left\{\frac{|A_l|^2}{ \langle N(f_l - f_p)\rangle^2}\right\} \geq \frac{\ln N}{N}$:
\begin{eqnarray}
r_L & \leq & \sqrt{L} O\left(\max_{ l,p; l \neq p} \left\{\frac{|A_l|}{ \langle N(f_l - f_p) \rangle}\right\}\right) \notag \\
& \leq & \sqrt{L} \max_{l} \left\{|A_l|\right\} O (N^{-1}) \notag
\end{eqnarray}
and the algorithm converges when $\max_{l} \left\{|A_l|\right\}/N$, which represents the maximum ratio of the amplitudes of the components and the signal length, is small enough so that $r_L<1$. It is important to emphasise that this result only holds when \textit{Assumption 1} holds the SNR is above the breakdown threshold of the algorithm. 
\end{enumerate}
\end{enumerate}
\end{thm}

%%%%%%%%%%%%%%%%%%%%%%%%%%%%%%%%%%%%%%
%
\begin{figure}[!t]
\centering{\includegraphics[width=3.45in]{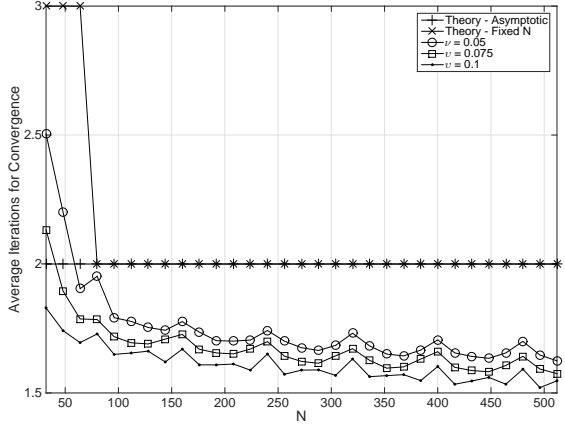}
\caption{Average iterations needed for the convergence of the proposed algorithm on (\ref{eq:sim_signal_1}). 5,000 Monte Carlo runs were used.}
\label{fig:converge}}
\end{figure}

\begin{figure}[!t]
\centering{\includegraphics[width=3.45in]{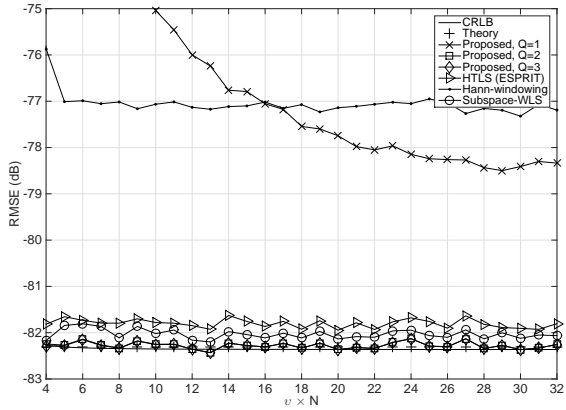}
\caption{RMSE of $\hat{f}_1$ versus $\upsilon N$  obtained by various methods on (\ref{eq:sim_signal_1}) when $\rho =20$dB and $a = 1$. 5,000 Monte Carlo runs were used.}
\label{fig:rmse_f1_theta}}
\end{figure}

\begin{figure}[!t]
\centering{\includegraphics[width=3.45in]{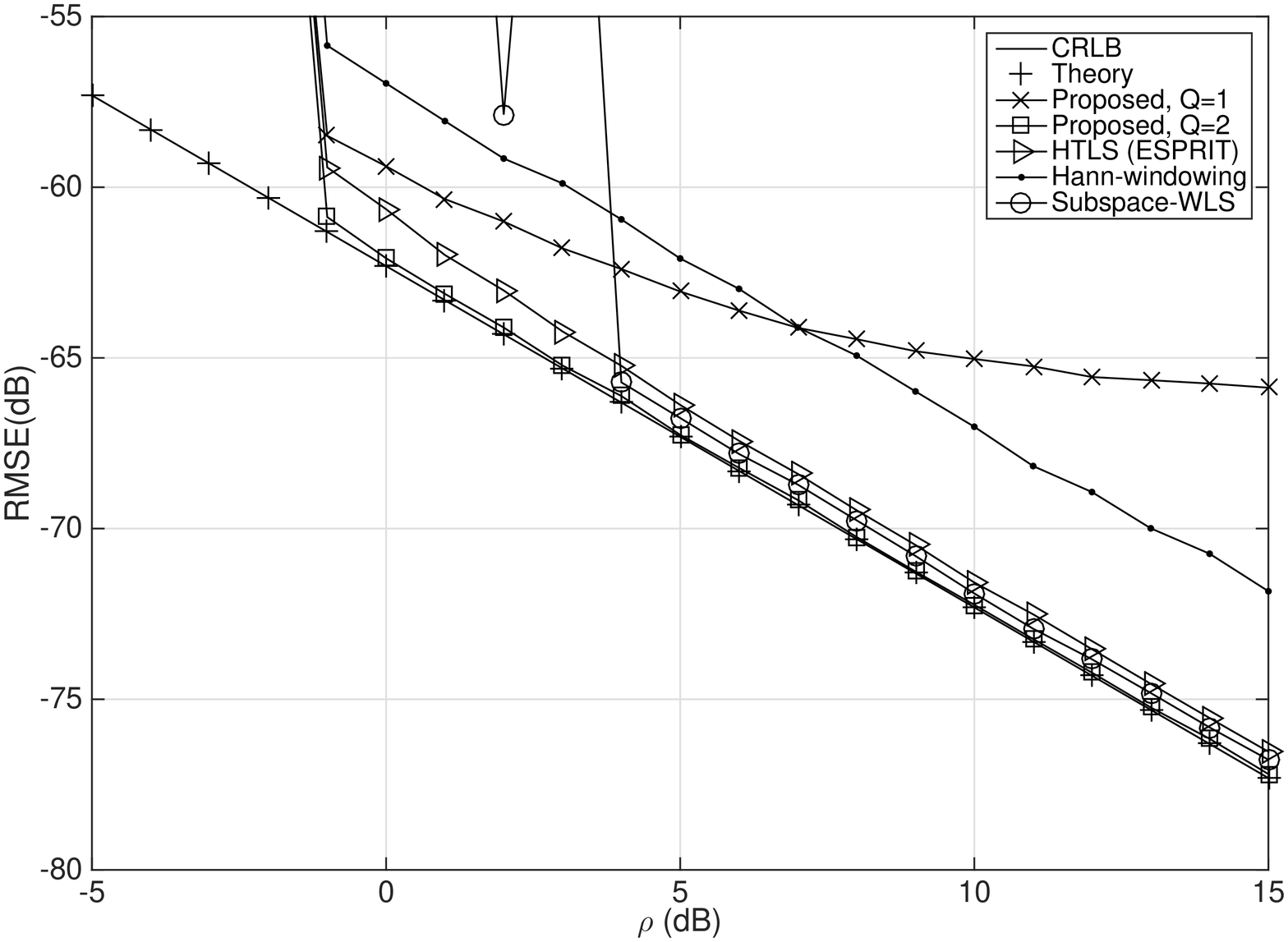}
\caption{RMSE of $\hat{f}_1$ versus $\rho$ using various estimation algorithms  on Eq. (\ref{eq:sim_signal_1}) when $\upsilon = 5/N$ and $a = 0.9$. 10,000 Monte Carlo runs were used.}
\label{fig:RMSE_f1}}
\end{figure}

\begin{figure}[!t]
\centering{\includegraphics[width=3.45in]{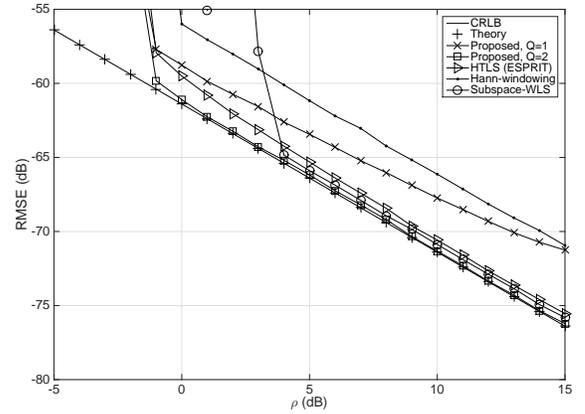} 
\caption{RMSE of $\hat{f}_2$ versus $\rho$ using various estimation algorithms on Eq. (\ref{eq:sim_signal_1}) when $\upsilon = 5/N$ and $a = 0.9$. 10,000 Monte Carlo runs were used.}
\label{fig:RMSE_f2}}
\end{figure}

\section{Simulation Results} \label{sec:Simu}
First we test the proposed algorithm on the following signal
\begin{equation} \label{eq:sim_signal_1}
x(n) = e^{ j2\pi f_1 n} + a e^{j\phi} e^{j2\pi (f_1 +\upsilon) + n} + w(n).
\end{equation}
Here $\upsilon$ is the interval between the two frequencies, with $\phi$ randomly selected in each run such that $\phi \in [-\pi, \pi]$.  We set $N = 64$. $a \leq 1$ is the ratio of the two magnitudes and the SNR of the first component is $\rho = 1 /\sigma^2$.

We start by verifying the convergence of the algorithm as we vary $N$ from $32$ to $512$. In this test we fix $a = 1$, and $f_1 = -0.48$ and randomly choose $\rho$ in each run in the range $0$ to $50$dB. Fig. \ref{fig:converge} shows the iterations needed for convergence as a function of $N$. we consider the convergence of the algorithm when the difference between the residual estimates of consecutively two iterations is less than the CRLB. In addition to the theoretical results, we given three curves of the average number of iterations needed for convergence for $\upsilon$ to be $0.05, 0.075$ and $0.01$. Keeping in mind that the theoretical results essentially give an upper bound on convergence, we see that the simulation curves are indeed bounded by the theoretical ones. 

Next we investigate the performance, in terms of the root mean square error (RMSE) versus $\upsilon$. In this test we set $\rho =20$dB, $a = 1$ and randomly select in each run $f \in [-0.5,0.5]$. We compare our algorithm with the Hankel Total Least Square (HTLS) method \cite{Vanhuffel1994}, frequency domain Hann-windowing method \cite{Diao2014} and subspace-WLS approach in \cite{Sun2012}. HTLS is essentially ESPRIT \cite{Haardt1995} with a Total Least Square (TLS) solution, and is the state-of-art time-domain parametric estimator. The Hann-windowing method and the subspace-WLS method are the most recently proposed approaches of their kinds respectively. We also include in Fig. \ref{fig:rmse_f1_theta} curves of the CRLB and theoretical variance. We show performance of our estimator for $Q = 1$, $2$. For HTLS, the degree of freedom is $M = [ N/3]$ for the best accuracy. For subspace-WLS, we set the rearranged matrix size to $16 \times 4$, which we found yields the best performance for randomly selected frequencies. First we see that $Q = 2$ is sufficient for our algorithm to reach CRLB-comparable performance at $ \upsilon N \geq 4$, and the resulting RMSE is approximately $1$dB less than that of HTLS. The Hann-windowing method exhibits the worst performance having an RMSE that is $5$dB higher than the CRLB. Although the case when $M$ is very small is beyond the scope of this paper, simulation results given in \cite{Ye2015} show that the estimator still exhibits excellent performance for $1.1\leq \upsilon N <4$.

% From Fig. \ref{fig:rmse_f1_theta2} we find when $a = 0.5$, that although the proposed method has higher RMSE than the high resolution methods at $\upsilon = 1/N$ due to slower convergence of $\delta_2$, it is capable of outperforming the other two methods at $1.1/N \leq \upsilon \leq 32/N$ using $Q \leq 25$.  

Now we show the RMSE versus $\rho$ for various frequency intervals. The simulation parameters are kept the same as the previous test. Figs. \ref{fig:RMSE_f1} and \ref{fig:RMSE_f2} give the RMSEs of $\hat{f}_1$ and $\hat{f}_2$ versus $\rho$ when $a = 0.9$ and $\upsilon = 5/N$. It is clear from the curves that the proposed estimator achieves the best performance in terms of RMSE and breakdown threshold after only $Q=2$ iterations.

%Then we investigate the performance versus $\rho$ using the following three-components signal
%%
%\begin{eqnarray} \label{eq:sim_signal_2}
%x(n) & = & e^{j2\pi f_1n} + 0.75e^{j\phi_2 + j2\pi (f_1 + \frac{5}{N})n} \notag \\
%&& \; + 0.5e^{j\phi_3 + j2\pi(f_2 + \frac{6}{N})n} + w(n), \quad n = 0,\ldots, 63.
%\end{eqnarray}
%%
%Here we set $a_1 = 1$, $f_2 = f_1 + 5/N$ and $f_3 = f_2 + 6/N$. $f_1$, $\phi_2$ and $\phi_3$ are randomly selected in each Monte Carlo run. Figs. \ref{fig:RMSE_f13}-\ref{fig:RMSE_f33} show the RMSEs of $\hat{f}_1$ to $\hat{f}_3$ respectively. These results show similar trends as the the two components case, confirming that our proposed method outperforms the HTLS, Hann-windowing and subspace-WLS algorithms. 

In order to show the performance advantage and robustness of the proposed algorithm we apply it to a signal with $L=15$ components:
\begin{equation} \label{eq:sim_signal_3}
x(n) =  \sum_{l = 1}^{15} A_l e^{j2\pi nf_l}  + w(n), \quad n = 0,\ldots, N-1,
\end{equation}
where $A_l = a_le^{j\phi_l}$, and $N=64$. The parameters $a_l$, $f_l$ (where $f_l = f_1 + \upsilon_l$) are generated randomly and fixed at the start of the simulation. We set without loss of generality $a_1 = 1$, and choose the other magnitudes randomly from the interval $0.25$ to $1$. The frequency separations $\upsilon_l$ are also uniformly distributed between $2/N$ and $4/N$ (which is essentially between $2$ and $4$ frequency bins). The SNR of the first component was set to $5$dB. The values of the various simulation parameters are given in table \ref{tab:SimulParam} along with the frequency separations in bins ($\upsilon N$) and SNRs. The results were obtained form $200,000$ Monte Carlo runs and in each run we generated the phases $\phi_l$ randomly from a uniform distribution over the interval $[-\pi,\pi]$.

\begin{table}[h!]
	\caption{Simulation Parameters for $L=15$ components.}
	\label{tab:SimulParam}
	\centering
	\begin{tabular}{|c|c|c|c|c|}
     \hline
     \hline
		Component& Amplitude&Frequency&Frequency&Nominal\\
		Number $l$&$a_l$&$f_l$&Separation $\upsilon N$ & SNR (dB)\\
		\hline
		1 & 1.0000 & -0.3071 & -& 5.0000\\
		2 & 0.6379 & -0.2623 & 2.8647 & 1.0947\\
		3 & 0.3825 & -0.2082 & 3.4616 & -3.3472\\
		4 & 0.8980 & -0.1609 & 3.0282 & 4.0651\\
		5 & 0.6046 & -0.1204 & 2.5947 & 0.6287\\
		6 & 0.9748 & -0.0855 & 2.2289 & 4.7784\\
		7 & 0.4310 & -0.0414 & 2.8284 & -2.3111\\
		8 & 0.5777 & -0.0080 & 2.1330 & 0.2344\\
		9 & 0.9284 & 0.0404 & 3.1010 & 4.3547\\
		10 & 0.8939 & 0.0785 & 2.4337 & 4.0262\\
		11 & 0.3282 & 0.1098 & 2.0082 & -4.6781\\
		12 & 0.4311 & 0.1655 & 3.5622 & -2.3080\\
		13 & 0.6182 & 0.2166 & 3.2689 & 0.8227\\
		14 & 0.8352 & 0.2683 & 3.3086 & 3.4360\\
		15 & 0.8690 & 0.3148 & 2.9780 & 3.7802\\
       \hline
       \hline
	\end{tabular}
\end{table}

%\textbf{I have removed the test on three components as it is no more informative than the two-component test. I will send you the results and the mfile to generate figures and I want you to put the following figures in:}
%
%\textbf{1. A figure of the entire spectrum: this figure shows that the algorithms exhibit concentrations around each of the 15 frequencies.}
%
%\textbf{2. Figures zooming in on some of the components chosen to show the difference in performance. I will look at the results and tell you which components.}
%
%\textbf{3. Figures showing the rmses of the algorithms which the mfile returns. I will then edit the English and finalise the paper. I would like it finished by the end of this week.}

\begin{figure}[!t]
\centering{\includegraphics[width=3.45in]{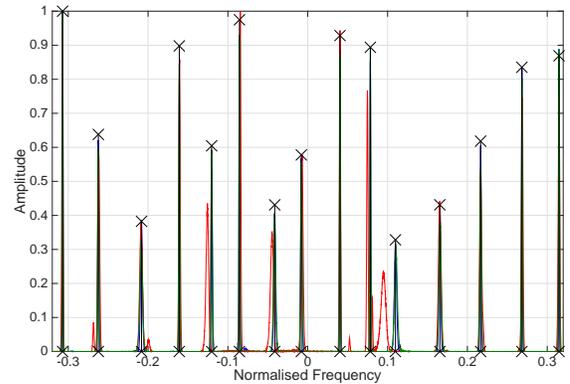} 
\caption{Normalised distributions of frequency estimates obtained by various methods on (\ref{eq:sim_signal_3}). 200,000 Monte Carlo runs were used. Blue: The proposed method using $Q = 3$; Red: The Hann-windowing method; Green: HTLS (ESPRIT); Black with crosses: True amplitudes at the true frequencies.}
\label{fig:multi_full}}
\end{figure}

\begin{figure}[htb]
\begin{minipage}[b]{0.48\linewidth}
  \centering
  \centerline{\includegraphics[width=1.7in]{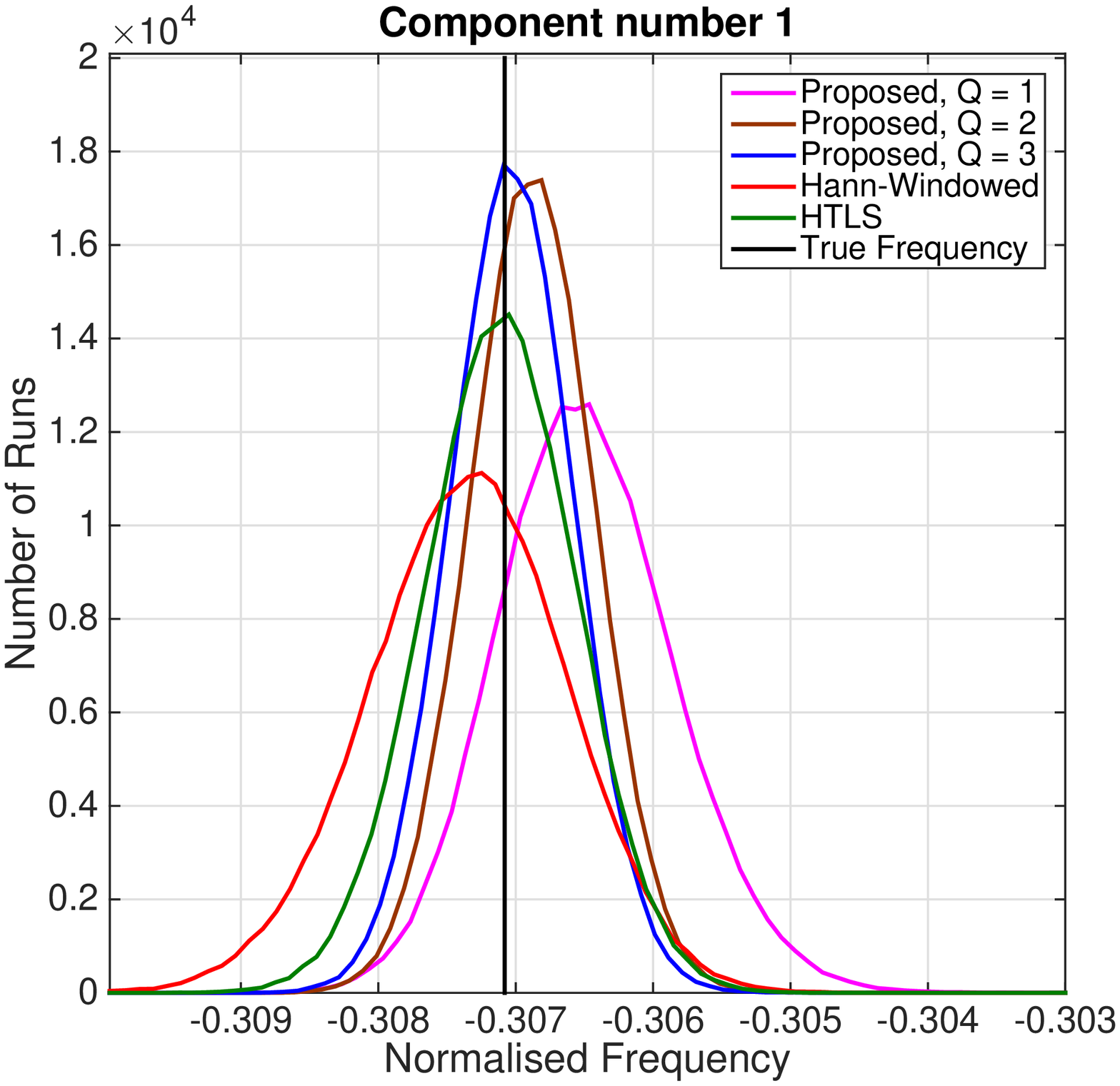}}
%  \vspace{2.0cm}
  \centerline{(a) Component No.1}\medskip
  \medskip
\end{minipage}
\hfill
\begin{minipage}[b]{0.48 \linewidth}
  \centering
  \centerline{\includegraphics[width=1.7in]{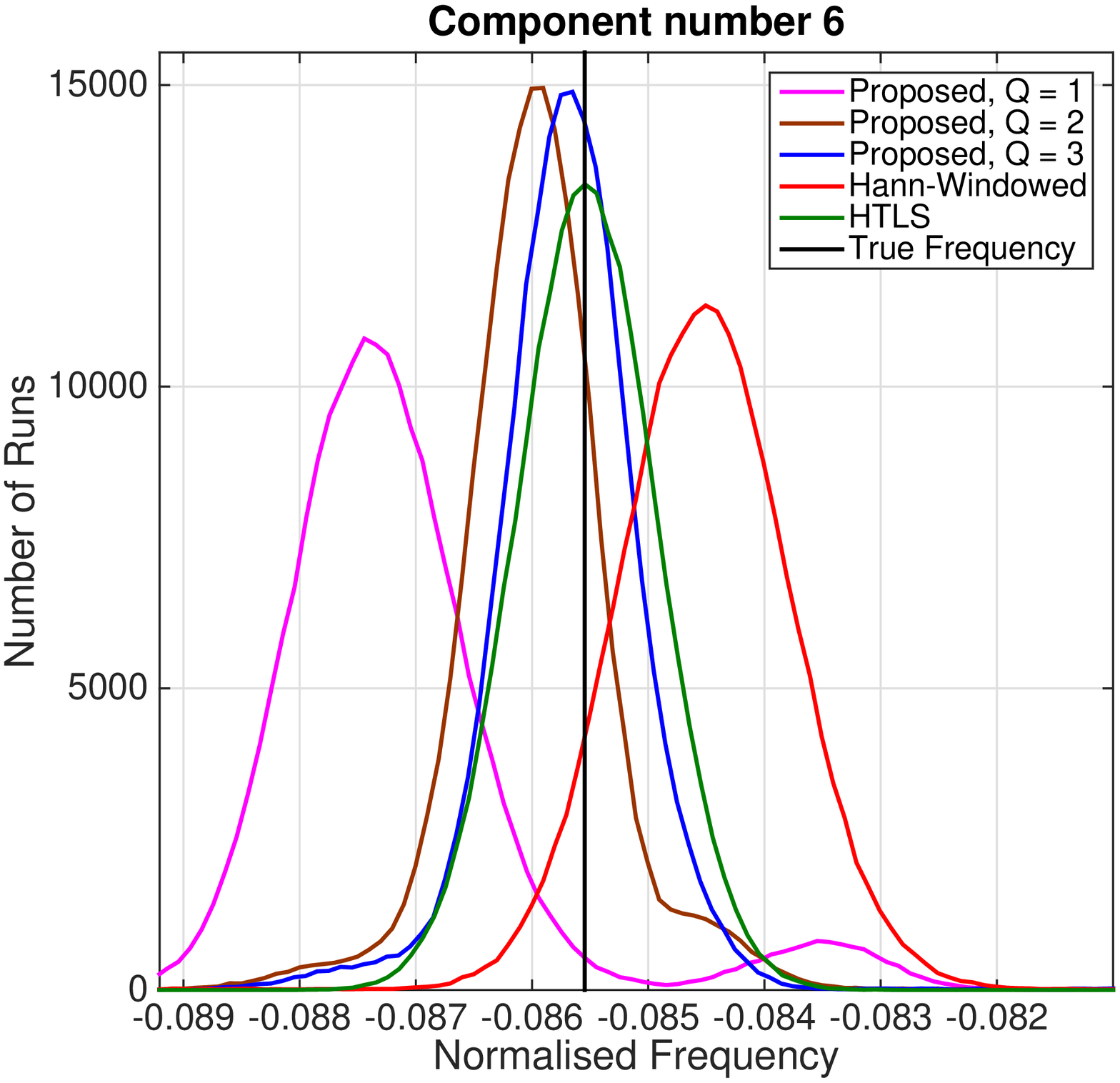}}
%  \vspace{2.0cm}
  \centerline{(b) Component No.6}\medskip
  \medskip
\end{minipage}
\begin{minipage}[b]{0.48 \linewidth}
  \centering
  \centerline{\includegraphics[width=1.7in]{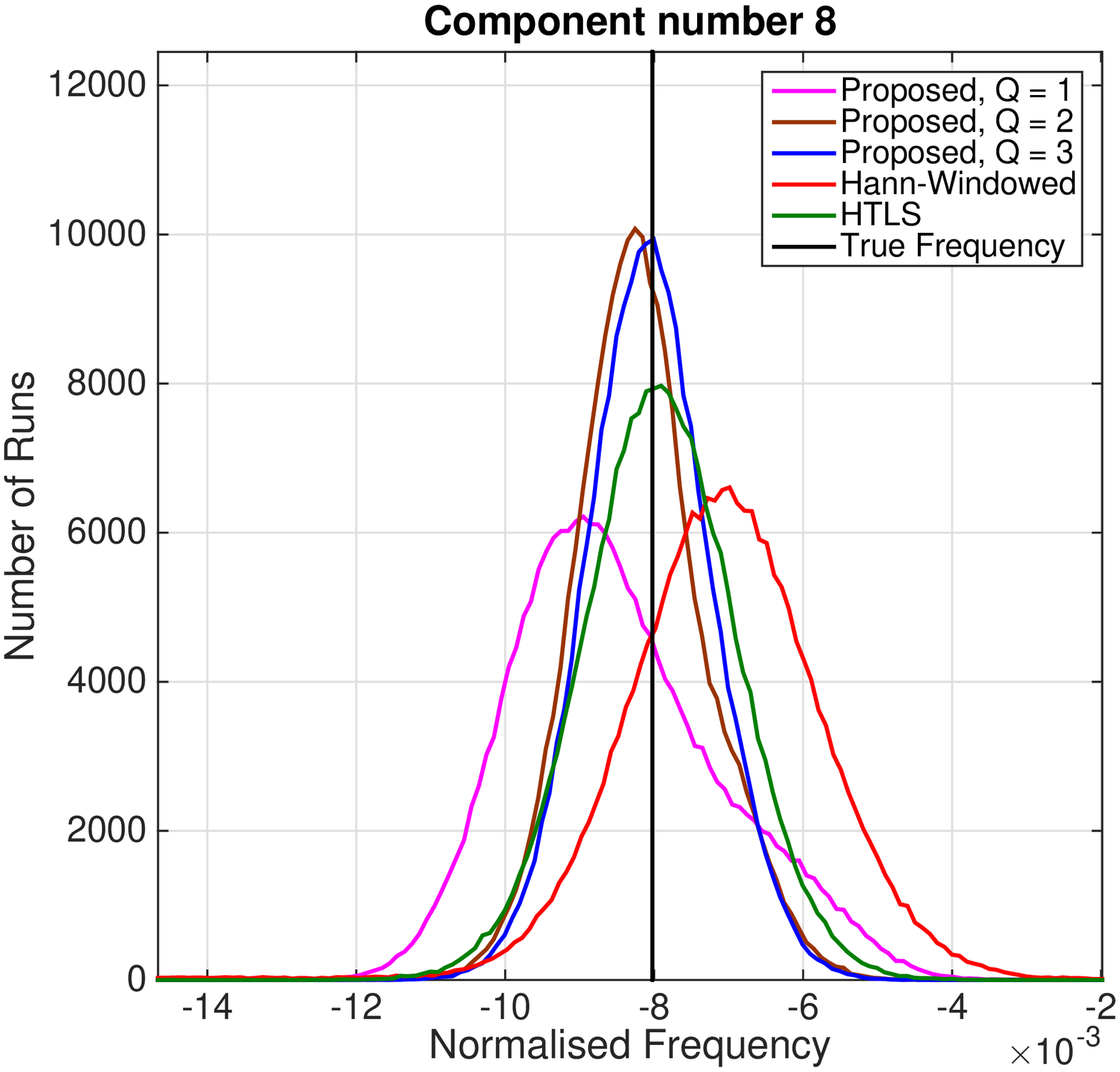}}
%  \vspace{2.0cm}
  \centerline{(c) Component No.8}\medskip
  \medskip
\end{minipage}
\hfill
\begin{minipage}[b]{0.48\linewidth}
  \centering
  \centerline{\includegraphics[width=1.7in]{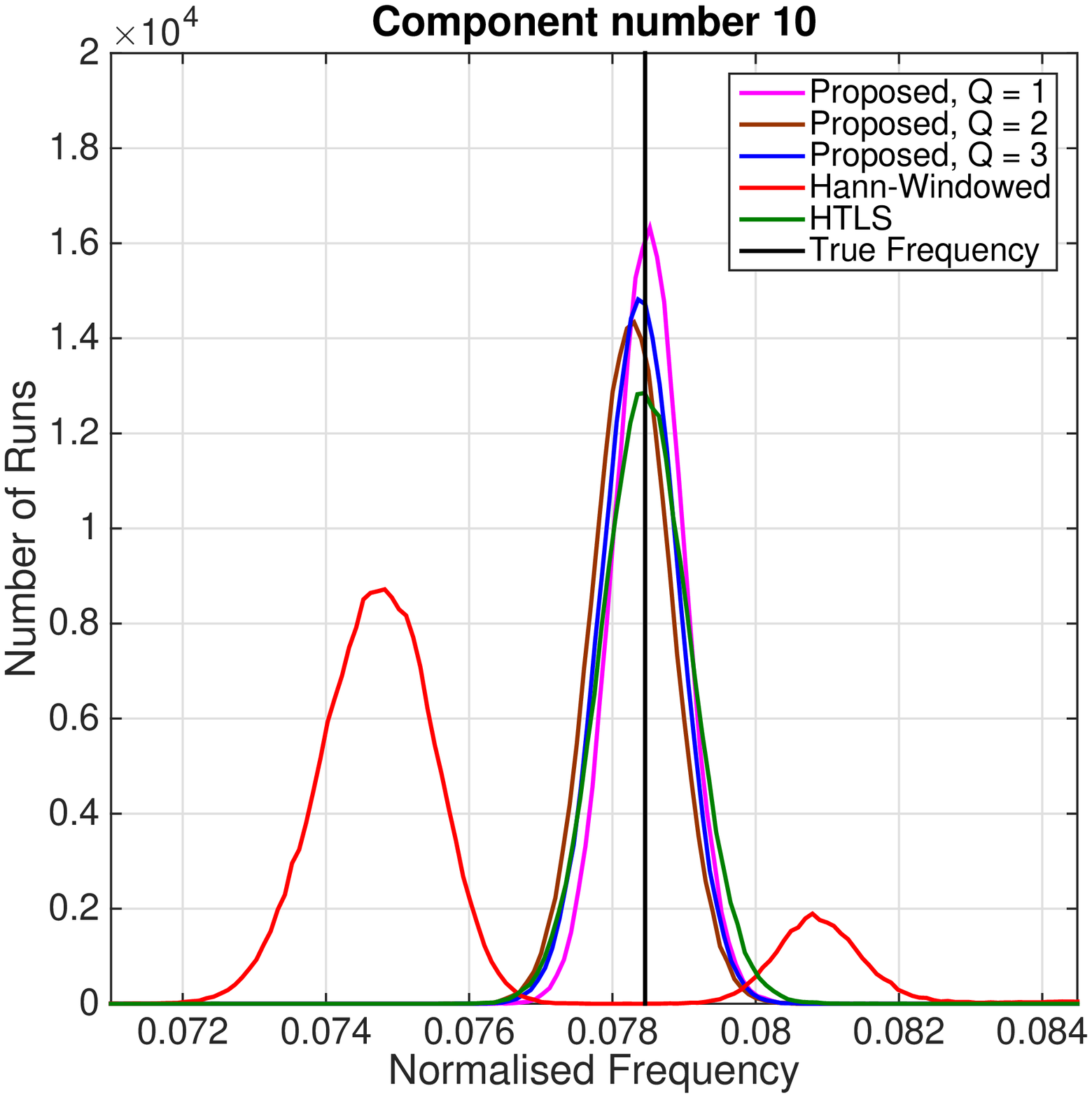}}
%  \vspace{2.0cm}
  \centerline{(d) Component No.10}\medskip
  \medskip
  \end{minipage}
\caption{Distributions of frequency estimates of varies components obtained by various methods on (\ref{eq:sim_signal_3}). 200,000 Monte Carlo runs were used.}
\label{fig:multi_comp}
\end{figure}

\begin{figure}[!t]
\centering{\includegraphics[width=3.45in]{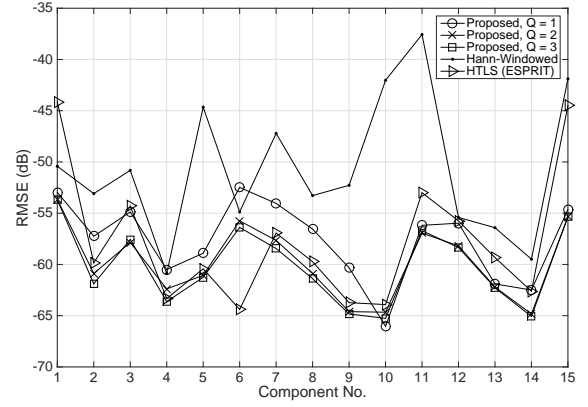} 
\caption{RMSE of frequency estimates obtained by various methods on (\ref{eq:sim_signal_3}). 200,000 Monte Carlo runs were used.}
\label{fig:multi_rmse}}
\end{figure}

\begin{figure}[!t]
\centering{\includegraphics[width=3.45in]{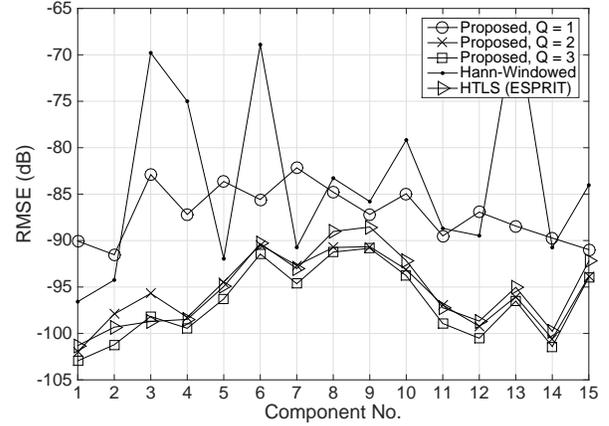} 
\caption{RMSE of frequency estimates obtained by various methods on (\ref{eq:sim_signal_3}) when $N = 1024$. 100,000 Monte Carlo runs were used.}
\label{fig:multi_rmse_1024}}
\end{figure}

Fig. \ref{fig:multi_full} shows the normalised distributions of the frequency estimates obtained by the proposed method, HTLS and the Hann-windowing method. The black markers represent the actual amplitudes at the true frequencies. Note that the subspace-WLS is not implementable in this test as the number of components is larger than the possible maximum length of the dimensions of the rearranged matrix. Although we ran the proposed method for $Q=6$ iterations, there was little change after 3 iterations. Therefore, we only show here the results up to $Q = 3$. We also set the degrees of freedom of HTLS to $M =  N/2$ which allows it to achieve the best performance. As we can see from the figure the proposed algorithm and HTLS achieve accurate estimation of the frequencies of all 15 components since the distributions of their estimates show peaks at each of the true frequencies. The Hann-windowed estimator, on the other hand, has the worst performance. Fig. \ref{fig:multi_comp} makes these observations clearer as we zoom-in on the distributions of components 1, 6, 8 and 10. The plots demonstrate that the estimates obtained by the proposed algorithm become more concentrated at the true frequencies as the number of iterations increases from 1 to 3. Fig. \ref{fig:multi_rmse} shows the RMSE of the frequency estimates of all the components obtained by the three methods. Again we see that the Hann-windowed method has the worst performance, whereas the proposed algorithm outperforms the HTLS algorithm at all but component 6. Similar results can be found in Fig. \ref{fig:multi_rmse_1024}, where we show the RMSE after increasing the signal length of Eq. (\ref{eq:sim_signal_3}) from $64$ to $1024$. 

Finally, we report the computational complexity of the estimators in table \ref{tab:ProcessingTimes}. Specifically, we give the ratio of the processing time of the HTLS algorithm to that of the proposed estimator. Although the Hann-Windowing method is very quick, it nonetheless does not match the estimation performance of the other two algorithms and so we ignore it. Despite our implementation of the proposed algorithm not being optimised, we see that as $N$ gets large, the gap between the two computational loads opens quite wide. For $N = 2048$, the HTLS algorithm takes almost 120 times longer than our algorithm to produce the estimates of the 15 components.

\begin{table}[h!]
	\caption{Ratio of the processing time $T_{HTLS}$ of the HTLS algorithm to that of the proposed estimator $T_{Multi}$.}
	\label{tab:ProcessingTimes}
	\centering
	\begin{tabular}{|l|l|l|l|l|}
     \hline
     \hline
		Number of Samples $N$& 256&512&1024&2048\\
		\hline
		$\frac{T_{HTLS}}{T_{Multi}}$&2.04&7.21&24.97&119.79\\
       \hline
       \hline
	\end{tabular}
\end{table} 

%\begin{figure}[!t]
%\centering{\includegraphics[width=4.5in]{r02}
%\caption{RMSE of $\hat{f}_2$ obtained by the proposed method versus $Q$ under various $\upsilon$ when SNR $=30$dB and $\alpha_2 = 1$. 1,000 Monte Carlo runs were used.}}
%\label{fig:rmse_f1_Q}
%\end{figure}

%\begin{figure}[!t]
%\centering{\includegraphics[width=4.5in]{r2}
%\caption{RMSE of $\hat{f}_2$ obtained by the proposed method versus $Q$ under various $\upsilon$ when SNR $=30$dB and $\alpha_2 = 0.5$. 1,000 Monte Carlo runs were used.}}
%\label{fig:rmse_f2_Q}
%\end{figure}

%\begin{figure}[!t]
%\centering{\includegraphics[width=3.45in]{QvsTheta}
%\caption{Numbers of iterations needed for convergence of the proposed algorithm versus frequency separation $\upsilon N$ when $a = 1$. 1,000 Monte Carlo runs were used.}}
%\label{fig:rmse_f1_Q}
%\end{figure}
%
%\begin{figure}[!t]
%\centering{\includegraphics[width=3.45in]{1_3}
%\caption{RMSE of $\hat{f}_1$ versus $\rho$ using various estimation algorithms on Eq. (\ref{eq:sim_signal_2}). 10,000 Monte Carlo runs were used.}
%\label{fig:RMSE_f13}}
%\end{figure}
%%
%\begin{figure}[!t]
%\centering{\includegraphics[width=3.45in]{2_3}
%\caption{RMSE of $\hat{f}_2$ versus $\rho$ using various estimation algorithms on Eq. (\ref{eq:sim_signal_2}). 10,000 Monte Carlo runs were used.}
%\label{fig:RMSE_f23}}
%\end{figure}
%%
%\begin{figure}[!t]
%\centering{\includegraphics[width=3.45in]{3_3}
%\caption{RMSE of $\hat{f}_3$ versus $\rho$ using various estimation algorithms on Eq. (\ref{eq:sim_signal_2}). 10,000 Monte Carlo runs were used.}
%\label{fig:RMSE_f33}}
%\end{figure}

\section{Conclusion}\label{sec:Conc}
In this paper, we proposed a novel method for accurately estimating the frequencies of multiple superimposed complex exponentials in noise. The proposed algorithm uses an efficient interpolation strategy to estimate the frequency of each component at a time in combination with an iterative leakage subtraction scheme. As the algorithm iterates, the leakage subtraction leads to a gradual reduction in the error of the interpolated coefficients due to the other components. Theoretical analysis showed that the algorithm converges to the asymptotic fixed point where the estimates are unbiased and minimum variance is just marginally larger than the CRLB. Simulation results were presented to verify the theoretical analysis. These results show that our method outperforms state-of-art methods in terms of estimation accuracy. Furthermore, the algorithm has a computational complexity that is same order  as the FFT operation, which is significantly lower than that of the non-parametric and time domain parametric estimators. 

\begin{appendices}
\section{Derivation of Theoretical Bias and Variance}
%We show in the following the derivation the theoretical bias and variance of $\hat{f}_1$ obtained by the proposed algorithm in the two component case. The related quantities of $\hat{f}_2$ can be easily obtained based on the presented results.
In the following derivation we consider the case $L=2$. Given $|f_1 - f_2| \sim O(1)$ we have that $M \sim O(N)$. Also $W_2 \sim O\left( N^{-\frac{1}{2}}\sqrt{\ln N}\right)$ a.s. \cite{An1983} and $\eta_\pm \sim O(N^{-1})$, which yields $W_{2\pm} \sim O\left( N^{-\frac{3}{2}}\sqrt{\ln N}\right)$ a.s.. We also have $W_{\pm} \sim O\left( N^{-\frac{1}{2}}\sqrt{\ln N}\right)$  a.s. and Eq. (\ref{eq:beta1}) implies that $\beta_{\pm} \sim O(1)$. Thus, the noise terms in $V_\pm$, $W_{2\pm}$ are of smaller order than $W_\pm$. Also, we know $V_\pm \sim O\left(N^{-1} \right)$ from Eq. (\ref{eq:V_final}). Now by putting:
\begin{eqnarray} \label{eq:Omega2}
\Omega  =  (V_+ + V_- + W_+ + W_-) (V_+ - V_- + W_+ - W_-),
\end{eqnarray}
the interpolation function $\hat{h}_1$ in Eq. (\ref{eq:h1}) becomes
\begin{eqnarray} 
\hat{h}_1 & = &\left\lbrace h_1 + j \frac{N\lambda_1}{4a_1z_1\sin\left(\frac{\pi}{N}\right)}(V_+ + V_- + W_{+}+W_{-})\right\rbrace \notag \\
&& \; \times \left\lbrace 1 -  j \frac{N\lambda_1}{2a_1z_1\sin\left(\frac{\pi}{N}\right)}(V_+ - V_-+W_{+}-W_{-}) \right\rbrace  + O(\Omega) \notag \\
& = & h_1 +  j\frac{N\lambda_1}{4a_1z_1\sin\left(\frac{\pi}{N}\right)} U + O(\Omega) \label{approx_h2},
\end{eqnarray}
where
\begin{eqnarray} \label{eq:U}
U = (1-2h_1) (V_+ + W_{+}) + (1+2h_1) (V_- + W_{-}).
\end{eqnarray}
Referring to Eq. (\ref{eq:gest}), the estimation error of $z^{-1}_1$ is
\begin{eqnarray}
\zeta_1 & = & \hat{z}_1^{-1} - z_1^{-1} \notag \\
& = & -2j(\hat{h}_1 - h_1)\sin\left(\frac{\pi}{N}\right) \notag \\
& = & \frac{N\lambda_1}{2a_1z_1} U + O\left(\frac{\Omega}{N}\right).
\end{eqnarray}
Expanding $\ln \hat{z}_1^{-1}$ as Taylor series around $\ln z_1^{-1}$ yields
\begin{eqnarray}
\ln \hat{z}_1^{-1} = -j\frac{2\pi}{N}\delta'_1 + z_1 \zeta_1 + O(z_1^2 \zeta^2_1),
\end{eqnarray} 
where $\delta'_1 = \delta_1 - \nu_1$. We then have
\begin{eqnarray} \label{eq:temp_delta1}
\hat{\delta}'_1 & = & -\frac{N}{2\pi} \Im \left\{ \ln \hat{z}_1^{-1} \right\} \notag \\
&  = &  \delta'_1 - \frac{N}{2\pi} \Im \left\{ z_1 \zeta_1 + O(z_1^2 \zeta^2_1) \right\}.
\end{eqnarray}
The estimation error of $\delta'_1$ can then be obtained by
\begin{eqnarray} \label{eq:varepsilon_1}
\varepsilon_1 & = & - \frac{N}{2\pi} \Im \left\{  z_1 \zeta_1 + O(z_1^2 \zeta_1^2) \right\} \notag  \\
& = & -\frac{N^2}{4\pi} \Im \left\{\frac{\lambda_1}{a_1}U\right\}+ O\left(\frac{\Omega}{N}\right).
\end{eqnarray}
To reach Eq. (\ref{eq:varepsilon_1}) we used the fact $z_1^2 \zeta_1^2 \sim O\left(\frac{\Omega}{N^2}\right)$. Setting $d_1 = \lambda_1e^{-j2\pi(\delta_1-\nu_1)}$, we have
\begin{eqnarray} \label{eq:varepsilon_1_fin}
\varepsilon_1 & = & -\frac{N^2}{8\pi A_1\cos(\pi\delta'_1)} \Im \left\{d_1U\right\}+ O\left(\frac{\Omega}{N}\right).
%& = & -\frac{N^2}{8\pi A_1\cos(\pi\delta'_1)}\left(d_{1R}U_I+d_{1I}U_R\right) \notag \\
%&& \;  + O\left(\frac{\Omega}{N}\right)\notag.
\end{eqnarray}
Asymptotically, the noise terms $W_{\pm}$ converge in distribution \cite{An1983}. Thus, $\hat{f}_1$ asymptotically converges to a normal:
\begin{equation}
\hat{f}_1 \rightarrow \mathcal{N}\left(f_1 + \mu_{f_1}, \sigma^2_{f_1}\right).
\end{equation} 

We know that $E[W_{\pm}] = E[W_{2\pm}] = 0$,  and the mean of (\ref{eq:Omega2}) only contains the contribution of the leakage terms $V_\pm$.  Now 
\begin{equation} \label{eq:h_approx}
h_1 = \delta'_1 + O(N^{-2}) = \delta_1 - \nu_1 + O(N^{-2}),
\end{equation}
then 
\begin{equation} \label{eq:lambda_approx}
\lambda_1  =  -\frac{4\pi^2}{N^2}\left[(\delta_1-\nu_1)^2 - 0.25 \right] + O(N^{-3}),
\end{equation}
and the mean of the asymptotic distribution of $\hat{f}_1$ is given by
\begin{eqnarray} \label{eq:Eep}
\mu_{{f}_1} & = & -\frac{N}{8\pi A_1 \cos(\pi\delta'_1)} \Im \left\{d_1\left[(1-2\delta'_1)\beta_+ + (1+2\delta'_1)\beta_-\right] \right\}\notag \\
& = &  \frac{\pi[(\delta_1 - \nu_1)^2 - 0.25]}{2 A_1 \cos(\pi(\delta_1-\nu_1))}\notag \\
&& \quad \left\{ [1-2(\delta_1-\nu_1)]\Im \left\{e^{-j2\pi(\delta_1 - \nu_1)}\beta_+\right\}\right. \notag \\
&& \qquad \left. + [1+2(\delta_1-\nu_1)]\Im \left\{e^{-j2\pi(\delta_1 - \nu_1)}\beta_-] \right\} \right\}.
\end{eqnarray}
Note that when $\delta_1 = \delta_2 = \nu_1 = \nu_2 = 0$ Eq. (\ref{eq:A2})  gives $\hat{A}_2 = A_2 + W_2$, which when  substituted into Eq. (\ref{eq:V_final}) yields $V_\pm = - W_{2\pm}$ and consequently $\mu_{{f}_1} |_{\delta_1 = \delta_2 = \nu_1 = \nu_2 = 0} = 0$. Now since
\begin{equation}
\textrm{Var}[W_2]  = \textrm{Var}[W_\pm] = \frac{\textrm{Var}[W_{2\pm}]}{|\eta_\pm|^2} =  \frac{\textrm{Var}[\eta_\pm W_{2}]}{|\eta_\pm|^2} = \frac{\sigma^2}{N} \notag,
\end{equation}
and the higher order term Eq. (\ref{eq:Omega2})  only contains the contribution of the noise terms, we arrive at
\begin{eqnarray}
\sigma^2_{f_1} &  = &  \frac{N^2}{64 \pi^2 A^2_1\cos^2(\pi\delta'_1)} \textrm{Var}[\Im\{d_1U\}].
\end{eqnarray}

Now turning to $\textrm{Var}[\Im\{d_1U\}]$, and using the fact $E[W_{2\pm}] = E[W_{\pm}] = 0$, we get
\begin{eqnarray} \label{eq:Var_temp0}
\textrm{Var} [\Im \{d_1 U\}] & = & \frac{1}{4} \textrm{Var} [d_1 U - d_1^* U^*] \notag \\
& = & \frac{1}{4} E[|(d_1U - E[d_1U])   - (d_1 U - E[d_1U])^*|^2] \notag \\
&=&  2|d_1|^2 E\left[|Z_+ + Z_-|^2\right] - d^2_1E[(Z_+ + Z_-)^2] \notag \\
&& \; - (d^*_1)^2E[(Z^*_+ + Z^*_-)^2]. 
\end{eqnarray}
where $Z_\pm  =  (1 \mp 2h_1) (W_\pm - W_{2\pm})$. Since$E[W^2_\pm] = E[W^2_2] = E[W_+ W_-] =  E[W_\pm W_2] = 0$, we have that $E[(Z_+ + Z_-)^2] = E[(Z^*_+ + Z^*_-)^2] = 0$ 
%\begin{eqnarray}
%E[(Z_\pm)^2] & = &   (1\mp 2h_1)^2 \left(E[W^2_\pm] + E[W^2_{2\pm}]  - 2E[W_{\pm}W_{2\pm}] \right) \notag \\
%& = &  (1\mp 2h_1)^2 \left(E[W^2_\pm] + E[W^2_{2\pm}]  - 2\eta_\pm E[W_{\pm}W_{2}] \right)  \notag \\
%& = & 0, \notag
%\end{eqnarray}
%and
%\begin{eqnarray}
%E[Z_+Z_-] & = & (1 - 4h^2_1) \left( E[(W_+ - W_{2+})(W_- - W_{2-})] \right) \notag \\
%& = & (1 - 4h^2_1) \left( E[W_+W_-] + \eta_+ \eta_- E[W^2_{2}] \right. \notag \\
%&& \; \left. - \eta_-E[W_{+}W_{2}] - \eta_+E[W_-W_{2}] \right) \notag \\
%& = & 0. \notag 
%\end{eqnarray}
%Now we focus on $E[(Z_+ + Z_-)^2]$
%Similarly we have 
%\begin{eqnarray}
%E[(W^*_\pm)^2] & = & E[(W^*_2)^2] = E[W^*_+ W^*_-] =  E[W^*_\pm W^*_2] \notag \\
%& = & 0 \notag,
%\end{eqnarray}
%and
%\begin{equation}
%E[(Z^*_+ + Z^*_-)^2] = 0.
%\end{equation}
and Eq.(\ref{eq:Var_temp0}) becomes
\begin{eqnarray}
\textrm{Var} [\Im \{d_1 U\}] & = & \frac{|\lambda_1|^2}{2}E\left[|Z_+ + Z_-|^2\right],
\end{eqnarray}
with $E\left[|Z_+ + Z_-|^2\right]$  given by Eq. (\ref{eq:EZ}).
\begin{figure*}
\begin{eqnarray} \label{eq:EZ}
E\left[|Z_+ + Z_-|^2\right] & = & |1-2h_1|^2 \textrm{Var}[W_{2+}] + |1+2h_1|^2 \textrm{Var}[W_{2-}] +|1-2h_1|^2 \textrm{Var}[W_{+}] + |1+2h_1|^2 \textrm{Var}[W_{-}] \notag \\
&& \; - |1-2h_1|^2 \left(E[W_{2+}W^*_{+}] + E[W^*_{2+} W_{+}] \right)  - |1+2h_1|^2 \left(E[W_{2-}W^*_{-}] + E[W^*_{2-} W_{-}] \right) \notag\\
&& \; +(1-2h_1)(1+2h^*_1) \left( E[W_+W^*_-] + E[W_{2+}W^*_{2-}] - E[W_{2+}W^*_-] - E[W^*_{2-}W_{+}] \right) \notag \\
&& \; +(1-2h^*_1)(1+2h_1) \left( E[W^*_+W_-] + E[W^*_{2+}W_{2-}] - E[W^*_{2+}W_-] - E[W_{2-}W^*_{+}] \right) .
\end{eqnarray}
\end{figure*}
%
%We know that
%\begin{equation}
%E[W_+W^*_-] = E[W^*_+W_-]  = 0, \notag
%\end{equation}
%
%\begin{eqnarray}
%E[W_{2+}W^*_{2-}]  =  \left(E[W^*_{2+}W_{2-}]\right)^* =  \eta_+\eta^*_- \frac{\sigma^2}{N}, \notag
%\end{eqnarray}
%\begin{eqnarray}
%E\left[W_2 W^*_\pm \right] & = & E\left[ \left(\frac{1}{N}\sum_k w(k) e^{-j2\pi \frac{m_2 + \nu_2}{N}} \right) \right. \notag \\
%&& \left.\qquad \left(\frac{1}{N}\sum_l w^*(l) e^{j2\pi\frac{m_1+ \nu_1 \pm 0.5}{N}} \right)\right] \notag \\
%& = &\frac{1}{N^2} \sum_k\sum_l E[w(k)w^*(l)] e^{-j2\pi\frac{M + \nu_2 - \nu_1 \mp 0.5}{N}}\\
%& = & E\left[ W^*_2 W_\pm  \right]  =   \frac{\sigma^2}{N}\eta^*_\pm, \notag 
%\end{eqnarray}
%\begin{eqnarray}
%E[W_{2+}W^*_+]  =  E[W^*_{2+}W_+] =   \eta_+ E[W^*_2 W_+]  =  |\eta_+|^2 \frac{\sigma^2}{N}, \notag 
%\end{eqnarray}
%\begin{eqnarray}
%E[W_{2-}W^*_-] = E[W^*_{2-}W_-] = |\eta_-|^2 \frac{\sigma^2}{N},\notag 
%\end{eqnarray}
%\begin{eqnarray}
%E[W_{2+}W^*_-]  = \left(E[W^*_{2+}W_-]\right)^* =  \eta_+ E[W_2 W^*_-] = \eta_+ \eta^*_- \frac{\sigma^2}{N}, \notag 
%\end{eqnarray}
%\begin{eqnarray}
%E[W_{2-}W^*_+]  =  \left(E[W^*_{2-}W_+]\right)^* = \eta_- E[W_2 W^*_+]  =  \eta^*_+ \eta_- \frac{\sigma^2}{N}. \notag 
%\end{eqnarray}

Thus, Eq. (\ref{eq:varepsilon_1}) leads to $\sigma^2_{f_1}$   being
\begin{eqnarray} \label{eq:Var_1}
\sigma^2_{f_1}& = & \frac{N^2}{64 \pi^2 A^2_1\cos^2(\pi\delta'_1)} \textrm{Var}[\Im\{d_1U\}] \notag \\
& = & \frac{\pi^2  [(\delta_1-\nu_1)^2-0.25]^2}{4\rho_1 N^3 \cos^2(\pi (\delta_1-\nu_1))} \left\{[1+4(\delta_1-\nu_1)^2] \right. \notag \\
&& \quad  - \frac{1}{2}[1-2(\delta_1-\nu_1)]^2 |\eta_+|^2 \notag \\
&& \quad - \frac{1}{2}[1+2(\delta_1-\nu_1)]^2 |\eta_-|^2 \notag \\
&& \quad \left. - \frac{1}{2}[1-4(\delta_1-\nu_1)^2] (\eta^*_+\eta_- + \eta_+ \eta^*_-)    \right\} .
\label{eq:var_fin}
\end{eqnarray}

Now putting $M = \alpha N$ where $\alpha <1$, $\eta_\pm$ becomes
\begin{eqnarray} 
\eta_\pm & = & \frac{1+e^{j2\pi (\nu_2-\nu_1)}}{N\left(1-e^{j\frac{2\pi}{N}(\alpha N+\nu_2 - \nu_1 \mp 0.5)}\right)} \notag \\
& = &  j\frac{1+e^{j2\pi (\nu_2-\nu_1)}}{2\pi \alpha N} + O(N^{-2})  \sim  O(M^{-1}). \label{eq:eta_fin}
\end{eqnarray}
As $M \sim O(N)$, the lower order terms involving $\eta_\pm$ can be ignored giving
\begin{equation}
\sigma^2_{f_1}  =  \frac{\pi^2  [(\delta_1-\nu_1)^2-0.25]^2}{4\rho_1 N^3\cos^2(\pi (\delta_1-\nu_1))} [1+4(\delta_1-\nu_1)^2].
\end{equation}

%\textbf{why the subscript 2? that is why $\alpha_2$ and not $\alpha$ which you use later? Also since $M = m_2-m_1$ and writing $m_p = [Nf_p]$, where $[\bullet]$ is rounding $\bullet$ to the nearest integer, then we have that $M = ([Nf_2]-[Nf_1])N$, implying that $\alpha = f_2-f_1$, which is the normalised frequency separation!} 

\section{Analysis of Convergence}
We proceed with to study the convergence of the proposed algorithm for the case of two components. Let 
\begin{eqnarray} \label{eq:gamma}
\gamma & = &  \frac{1}{N}\frac{1 - e^{j2\pi(\delta_2 - \nu_2)}}{1 - e^{j\frac{2\pi}{N}(\delta_2 - \nu_2)}} \notag \\
& = & \frac{1}{N}\frac{1 - e^{j2\pi(\delta_2 - \nu_2)}}{-j\frac{2\pi}{N}(\delta_2 - \nu_2) (1 - O(N^{-1}))} \notag \\
& = & \frac{1 - e^{j2\pi(\delta_2 - \nu_2)}}{-j2\pi(\delta_2 - \nu_2)} + O(N^{-1}) \notag \\
& = & 1 + \delta_2 - \nu_2 + O((\delta_2 - \nu_2)^2),
\end{eqnarray}
where in these manipulations weused the Taylor expansions $1 - e^x = - x + O(x^2)$ and $1/(1+x) = 1-x + O(x^2)$. Using the same expansions and putting $M = \alpha N$, $\beta_\pm$ becomes
\begin{eqnarray}
\beta_\pm &  = & \frac{1+e^{j2\pi(\delta_2 - \nu_1)}}{1-e^{\frac{j2\pi}{N}M} (1 + O(N^{-1}))}  - \gamma  \frac{1+e^{j2\pi(\nu_2 - \nu_1)}}{1-e^{\frac{j2\pi}{N}M} (1 + O(N^{-1}))} \notag \\
& = &  \frac{1+e^{j2\pi(\delta_2 - \nu_1)} - \gamma [1+e^{j2\pi(\nu_2 - \nu_1)}]}{1-e^{\frac{j2\pi}{N}M}} + O(N^{-1})  \notag \\
& = & j\frac{1}{2\pi \left( 1-e^{\frac{j2\pi}{N}M} \right)}  \left\{(\delta_2 - \nu_1) + O((\delta_2 - \nu_1)^2) \right. \notag \\
&& \quad  - [1+(\delta_2 - \nu_2) + O((\delta_2 - \nu_2)^2)] \notag \\
&& \qquad \left. \times [(\nu_2 - \nu_1)+ O((\nu_2 - \nu_1)^2)]\right\} + O(N^{-1})  \notag \\
 & = &  \frac{A_2}{\alpha} (\delta_2 - \nu_2) O(1), 
\end{eqnarray}
so we have 
\begin{eqnarray}
V_\pm  =  \frac{A_2}{N}\beta_\pm + W_{2\pm} + O(N^{-2}) = \frac{A_2}{\alpha}(\delta_2 - \nu_2) O(N^{-1}).
\end{eqnarray}
Now we rewrite $U$ in Eq. (\ref{eq:U}) as
\begin{eqnarray}
U & = & (1-2h_1) V_+ + (1+2h_1) V_- \notag \\
&& \; + (1-2h_1)  W_{+} + (1+2h_1) W_{-} \notag \\
& = & U_V + U_W,
\end{eqnarray}
and $\varepsilon_1$ in Eq. (\ref{eq:varepsilon_1}) becomes
\begin{eqnarray}
\varepsilon_1 = \varepsilon_{1,V} + \varepsilon_{1,W} +  O(\Omega).
\end{eqnarray}
It is clear that $\varepsilon_{1,W}$ is identical to the estimation error of the single component case. From \cite{Aboutanios2005} and \cite{Aboutanios2002} we know that
\begin{equation}
\varepsilon = \left[ \delta_1 +  (\nu_1 - \delta_1) O(N^{-\frac{1}{2}}\sqrt{\ln N})\right] O(N^{-\frac{1}{2}}\sqrt{\ln N}). \quad \textrm{a.s.}
\end{equation}
On the other hand, we have
\begin{eqnarray}
\varepsilon_{1,V}  = \frac{A_2}{\alpha}(\delta_2 - \nu_2) O(N^{-1}).
\end{eqnarray}

Now the iterative estimation function of $\delta_1$ can be written as a function of $\nu_1, \nu_2$ as
\begin{eqnarray}
H_1(\nu_1, \nu_2) = \nu_1 + (\delta_1 - \nu_1) + \varepsilon_1(\nu_1, \nu_2) = \delta_1 + \varepsilon_1(\nu_1, \nu_2).
\end{eqnarray}
For any $\nu_1^{(1)}, \nu_2^{(1)}, \nu_1^{(2)}, \nu_2^{(2)} \in [-0.5, 0.5]$ we have
\begin{eqnarray}
|H_1(\nu_1^{(2)}, \nu_2^{(2)}) & - & H_1(\nu_1^{(1)}, \nu_2^{(1)})| \notag \\
& \leq &  \frac{A_2}{\alpha}|\nu_2^{(2)} - \nu_2^{(1)}| O(N^{-1}) \notag \\
&& \; + |\nu_1^{(2)}  - \nu_1^{(1)}| O( N^{-\frac{1}{2}}\sqrt{\ln N}) \notag \\
&& \; + \left( |\nu_2^{(2)} - \nu_2^{(1)}| +  |\nu_1^{(2)} - \nu_1^{(1)}| \right)  O(\Omega) \notag \\
& \leq & \left( |\nu_2^{(2)} - \nu_2^{(1)}| +  |\nu_1^{(2)} - \nu_1^{(1)}| \right) O(\sqrt{\Lambda_{2}}),  \quad \textrm{a.s.} \notag
\end{eqnarray}
with $\Lambda_2$ given by
\begin{eqnarray} \label{eq:Lambda}
\Lambda_2 = \max \left \{ \frac{|A_2|^2}{\langle(f_2 - f_1)N\rangle^2}, \frac{\ln N}{N}\right\}.
\end{eqnarray} 
To reach Eq. (\ref{eq:Lambda}) we used the facts that $\rho_l = |A_l|^2/\sigma^2$ and 
\begin{eqnarray}
f_2 - f_1 = \frac{\alpha N + \delta_2 - \delta_1}{N} & \Leftrightarrow & \alpha N = \langle N(f_2 - f_1)\rangle.
\end{eqnarray}
Proceeding similarly for $\delta_2$ yields 
\begin{eqnarray}
 |H_2(\nu_1^{(2)}, \nu_2^{(2)}) & - & H_2(\nu_1^{(1)}, \nu_2^{(1)})| \notag \\
& \leq & \left( |\nu_2^{(2)} - \nu_2^{(1)}| +  |\nu_1^{(2)} - \nu_1^{(1)}| \right) O(\sqrt{\Lambda_1}),  \quad \textrm{a.s.} \notag 
\end{eqnarray}
with $\Lambda_1 =  \max \left \{ \frac{|A_1|^2}{\langle(f_1 - f_2)N\rangle^2}, \frac{\ln N}{N}\right\}$. Finally we have
\begin{eqnarray} \label{eq:contractive}
&& |H_1(\nu_1^{(2)}, \nu_2^{(2)})  - H_1(\nu_1^{(1)}, \nu_2^{(1)})|^2 \notag \\
&& \quad + |H_2(\nu_1^{(2)}, \nu_2^{(2)}) - H_2(\nu_1^{(1)}, \nu_2^{(1)})|^2 \notag \\
& \leq &  2 (O(\Lambda_1) + O(\Lambda_2)) \left( |\nu_2^{(2)} - \nu_2^{(1)}|^2 +  |\nu_1^{(2)} - \nu_1^{(1)}|^2 \right) \notag \\
& \leq &  2 O\left(\max\left\{\Gamma_{1,2}, \Gamma_{2,1}, \frac{\ln N}{N} \right\}\right) \notag \\
&& \qquad \left( |\nu_2^{(2)} - \nu_2^{(1)}|^2 +  |\nu_1^{(2)} - \nu_1^{(1)}|^2 \right), \quad \textrm{a.s.}
\end{eqnarray}
Here $\Gamma_{l,p}$ is given by Eq. (\ref{eq:Gamma2}). 
Eq. (\ref{eq:contractive}) implies that the estimation procedure is a contractive mapping and according to the fixed point theorem \cite{Aboutanios2005}, the iterative estimator converges asymptotically to the fixed point of the true frequency residual:
\begin{equation} \label{eq:gen_contractive}
(H_1(\delta_1,\delta_2), H_2(\delta_1,\delta_2)) = (\delta_1, \delta_2),
\end{equation}
with the upper bound on the convergence rate being
\begin{equation}
r_2 \leq \sqrt{2}O\left( \max\left\{\sqrt{\Gamma_{1,2}}, \sqrt{\Gamma_{2,1}}, \sqrt{\frac{\ln N}{N}}\right\}\right).  \quad \textrm{a.s.}
\end{equation}
Consequently, the asymptotic mean and variance of $\hat{f}_1$ in the two-component case can be obtained by substituting $\delta_1 = \delta_2 = \nu_1 = \nu_2 = 0$ in Eqs. (\ref{eq:E}) and (\ref{eq:Var}), which results in
\begin{equation} \label{eq:varf2}
\mu_{f_1} = 0, \quad \textrm{and} \quad \sigma^2_{f_1} =  \frac{\pi^2}{64\rho_1 N^3}.
\end{equation}
Thus, the algorithm is asymptotically unbiased and achieves the minimum variance of the single component case. 

The above argument for convergence can be extended to the general case when $L\geq 2$ giving
\begin{eqnarray} 
&& \sum_{l = 1}^{L}|H_l(\nu_1^{(2)}, \ldots, \nu_L^{(2)})  - H_l(\nu_1^{(1)}, \ldots, \nu_L^{(1)})|^2 \notag \\
%& \leq &  2\sum_{l=1}^L O(\Gamma_l^2) \sum_{l =1}^L |\nu_l^{(2)} - \nu_l^{(1)}|^2\notag \\
& \leq &  L O\left( \max_{l,p = 1\ldots L; l \neq p} \left\{\Gamma_{l,p}, \frac{\ln N}{N} \right\}\right) \sum_{l =1}^L |\nu_l^{(2)} - \nu_l^{(1)}|^2. \; \textrm{a.s.}
\end{eqnarray}
The algorithm then converges to the $L$-dimensional fixed point
\begin{equation}
(H_1(\delta_1\ldots \delta_L), \ldots, H_L(\delta_1\ldots \delta_L)) = (\delta_1 \ldots \delta_L),
\end{equation}
and the upper bound of the convergence rate is
\begin{equation}
r_L \leq \sqrt{L}O\left( \max_{l,p = 1\ldots L; l \neq p}\left\{\sqrt{\Gamma_{l,p}}, \sqrt{\frac{\ln N}{N}}\right\}\right). \quad \textrm{a.s.}
\end{equation}
Furthermore, the mean and variance of the asymptotic distribution of $\hat{f}_p$ are obtained by substituting $\delta_l = \nu_l = 0, l = 1\ldots L$ in Eqs. (\ref{eq:E_gen}) and (\ref{eq:var_gen}) giving
\begin{eqnarray} 
\mu_{f_p} = 0, \quad \textrm{and} \quad \sigma^2_{f_p} =  \frac{\pi^2}{64\rho_p N^3}.
\end{eqnarray}

Now given \textit{Assumption 1}, the ACRLB for multiple component case equals to the single component case is, \cite{Yao1995},
\begin{equation}
\textrm{ACRLB}_{f_p} = \frac{6}{4\pi^2\rho_p N^3},
\end{equation}
Then the asymptotic ratio of $ \sigma^2_{f_p}$ to $\textrm{ACRLB}_{f_p}$ is given by $R = (0.25\pi^2)^2/6 \approx 1.0147$.

\end{appendices}
\bibliographystyle{IEEEtran}
\bibliography{refs}

\begin{thebibliography}{10}
\providecommand{\url}[1]{#1}
\csname url@samestyle\endcsname
\providecommand{\newblock}{\relax}
\providecommand{\bibinfo}[2]{#2}
\providecommand{\BIBentrySTDinterwordspacing}{\spaceskip=0pt\relax}
\providecommand{\BIBentryALTinterwordstretchfactor}{4}
\providecommand{\BIBentryALTinterwordspacing}{\spaceskip=\fontdimen2\font plus
\BIBentryALTinterwordstretchfactor\fontdimen3\font minus
  \fontdimen4\font\relax}
\providecommand{\BIBforeignlanguage}[2]{{%
\expandafter\ifx\csname l@#1\endcsname\relax
\typeout{** WARNING: IEEEtran.bst: No hyphenation pattern has been}%
\typeout{** loaded for the language `#1'. Using the pattern for}%
\typeout{** the default language instead.}%
\else
\language=\csname l@#1\endcsname
\fi
#2}}
\providecommand{\BIBdecl}{\relax}
\BIBdecl

\bibitem{Duda2011}
K.~Duda, L.~Magalas, M.~Majewski, and T.~Zielinski, ``{DFT}-based estimation of
  damped oscillation parameters in low-frequency mechanical spectroscopy,''
  \emph{IEEE Transactions on Instrumentation and Measurement}, vol.~60, no.~11,
  pp. 3608 --3618, 2011.

\bibitem{Umesh1996}
S.~Umesh and D.~Tufts, ``Estimation of parameters of exponentially damped
  sinusoids using fast maximum likelihood estimation with application to {NMR}
  spectroscopy data,'' \emph{IEEE Transactions on Signal Processing}, vol.~44,
  no.~9, pp. 2245--2259, 1996.

\bibitem{Stoica1997}
P.~Stoica and R.~L. Moses, \emph{Introduction to spectral analysis}.\hskip 1em
  plus 0.5em minus 0.4em\relax Prentice hall Upper Saddle River, 1997.

\bibitem{Zielinski2011}
T.~Zielinski and K.~Duda, ``Frequency and damping estimation methods - an
  overview,'' \emph{Metrology and Measurement Systems}, vol.~18, no.~4, pp.
  505--528, 2011.

\bibitem{Capon1969}
J.~Capon, ``High-resolution frequency-wavenumber spectrum analysis,''
  \emph{Proceedings of the IEEE}, vol.~57, no.~8, pp. 1408 -- 1418, 1969.

\bibitem{Li1996}
J.~Li and P.~Stoica, ``An adaptive filtering approach to spectral estimation
  and {SAR} imaging,'' \emph{IEEE Transactions on Signal Processing}, vol.~44,
  no.~6, pp. 1469--1484, 1996.

\bibitem{Yardibi2010}
T.~Yardibi, J.~Li, P.~Stoica, M.~Xue, and A.~Baggeroer, ``Source localization
  and sensing: {A} nonparametric iterative adaptive approach based on weighted
  least squares,'' \emph{IEEE Transactions on Aerospace and Electronic
  Systems}, vol.~46, no.~1, pp. 425--443, Jan 2010.

\bibitem{Glentis2008}
G.-O. Glentis, ``A fast algorithm for \textnormal{APES} and \textnormal{Capon}
  spectral estimation,'' \emph{IEEE Transactions on Signal Processing},
  vol.~56, no.~9, pp. 4207--4220, 2008.

\bibitem{Angel2012}
K.~Angelopoulos, G.~Glentis, and A.~Jakobsson, ``Computationally efficient
  {Capon-} and {APES-based} coherence spectrum estimation,'' \emph{IEEE
  Transactions on Signal Processing}, vol.~60, no.~12, pp. 6674--6681, Dec
  2012.

\bibitem{Ming2011}
M.~Xue, L.~Xu, and J.~Li, ``{IAA} spectral estimation: Fast implementation
  using the gohberg-semencul factorization,'' \emph{IEEE Transactions on Signal
  Processing}, vol.~59, no.~7, pp. 3251--3261, July 2011.

\bibitem{Hua1990}
Y.~Hua and T.~K. Sarkar, ``Matrix pencil method for estimating parameters of
  exponentially damped/undamped sinusoids in noise,'' \emph{IEEE Transactions
  on Acoustics, Speech, and Signal Processing}, vol.~38, no.~5, pp. 814--824,
  1990.

\bibitem{Haardt1995}
M.~Haardt and J.~A. Nossek, ``Unitary {ESPRIT}: How to obtain increased
  estimation accuracy with a reduced computational burden,'' \emph{IEEE
  Transactions on Signal Processing}, vol.~43, no.~5, pp. 1232--1242, 1995.

\bibitem{Vanhuffel1994}
S.~Vanhuffel, H.~Chen, C.~Decanniere, and P.~Vanhecke, ``Algorithm for
  time-domain \textnormal{NMR} data fitting based on total least squares,''
  \emph{Journal of Magnetic Resonance, Series A}, vol. 110, no.~2, pp. 228 --
  237, 1994.

\bibitem{Bresler1986}
Y.~Bresler and A.~Macovski, ``Exact maximum likelihood parameter estimation of
  superimposed exponential signals in noise,'' \emph{IEEE Transactions on
  Acoustics, Speech and Signal Processing}, vol.~34, no.~5, pp. 1081--1089, Oct
  1986.

\bibitem{Chan2010}
F.~Chan, H.~So, W.~Lau, and C.~Chan, ``Efficient approach for sinusoidal
  frequency estimation of gapped data,'' \emph{IEEE Signal Processing Letters},
  vol.~17, no.~6, pp. 611--614, June 2010.

\bibitem{Sun2012}
W.~Sun and H.~So, ``Efficient parameter estimation of multiple damped sinusoids
  by combining subspace and weighted least squares techniques,'' in
  \emph{ICASSP, IEEE International Conference on Acoustics, Speech and Signal
  Processing - Proceedings}, 2012, pp. 3509--3512.

\bibitem{Gough1994}
P.~T. Gough, ``A fast spectral estimation algorithm based on the {FFT},''
  \emph{IEEE Transactions on Signal Processing}, vol.~42, no.~6, pp.
  1317--1322, 1994.

\bibitem{Aboutanios2011}
E.~Aboutanios, ``Estimating the parameters of sinusoids and decaying sinusoids
  in noise,'' \emph{IEEE Instrumentation and Measurement Magazine}, vol.~14,
  no.~2, pp. 8--14, 2011.

\bibitem{Aboutanios2005}
E.~Aboutanios and B.~Mulgrew, ``Iterative frequency estimation by interpolation
  on {Fourier} coefficients,'' \emph{IEEE Transactions on Signal Processing},
  vol.~53, no.~4, pp. 1237--1242, 2005.

\bibitem{Quinn1994}
B.~G. Quinn, ``Estimating frequency by interpolation using {Fourier}
  coefficients,'' \emph{IEEE Transactions on Signal Processing}, vol.~42,
  no.~5, pp. 1264--1268, 1994.

\bibitem{Yang2011}
C.~Yang and G.~Wei, ``\BIBforeignlanguage{English}{A noniterative frequency
  estimator with rational combination of three spectrum lines},''
  \emph{\BIBforeignlanguage{English}{IEEE Transactions on Signal Processing}},
  vol.~59, no.~10, pp. 5065--5070, 2011.

\bibitem{Duda2014}
K.~Duda and S.~Barczentewicz, ``Interpolated {DFT} for $sin^{\alpha}(x)$
  windows,'' \emph{IEEE Transactions on Instrumentation and Measurement},
  vol.~63, no.~4, pp. 754--760, 2014.

\bibitem{Diao2014}
R.~Diao and Q.~Meng, ``An interpolation algorithm for discrete {Fourier}
  transforms of weighted damped sinusoidal signals,'' \emph{IEEE Transactions
  on Instrumentation and Measurement}, vol.~63, no.~6, pp. 1505--1513, 2014.

\bibitem{Quinn2001}
B.~G. Quinn and E.~J. Hannan, \emph{{T}he estimation and tracking of
  frequency}.\hskip 1em plus 0.5em minus 0.4em\relax New York: Cambridge
  University Press, 2001.

\bibitem{Belega2014}
D.~Belega and D.~Petri, ``Sine-wave parameter estimation by interpolated {DFT}
  method based on new cosine windows with high interference rejection
  capability,'' \emph{Digital Signal Processing}, vol.~33, no.~0, pp. 60 -- 70,
  2014.

\bibitem{An1983}
H.-Z. An, Z.-G. Chen, and E.~Hannan, ``The maximum of the periodogram,''
  \emph{Journal of Multivariate Analysis}, vol.~13, no.~3, pp. 383 -- 400,
  1983.

\bibitem{Yao1995}
Y.~Yao and S.~M. Pandit, ``{Cramer-Rao} lower bounds for a damped sinusoidal
  process,'' \emph{IEEE Transactions on Signal Processing}, vol.~43, no.~4, pp.
  878--885, 1995.

\bibitem{Ye2015}
S.~Ye and E.~Aboutanios, ``An algorithm for the parameter estimation of
  multiple superimposed exponentials in noise,'' in \emph{ICASSP, The
  International Conference on Acoustic, Speech and Signal Processing -
  Proceedings}, 2015.

\bibitem{Aboutanios2002}
E.~Aboutanios, ``Frequency estimation for low earth orbit satellites,'' Ph.D.
  dissertation, Univ. Technology, Sydney, Australia, 2002.

\end{thebibliography}

\end{document}